\documentclass[11pt,a4paper,reqno]{amsart}
\usepackage{color,latexsym,amsfonts,amssymb,bbm,comment}
\usepackage{hyperref}
\usepackage{amsmath,cite}
\usepackage{amsthm}
\usepackage{fullpage}
\usepackage{graphicx}
\usepackage{dsfont}
\usepackage{caption,subcaption}

\newtheorem {thm}{Theorem}[section]
\newtheorem {prop}[thm]{Proposition} 

\newtheorem {lem}[thm]{Lemma}

\newtheorem {defn}[thm]{Definition}

\newtheorem {cond}[thm]{Condition}

\newenvironment{remark}[1][Remark:]{\begin{trivlist}
\item[\hskip \labelsep {\bfseries #1}]}{\end{trivlist}}

\def\N{{\Bbb N}}

\def\Z{{\Bbb Z}}

\def\one{{\mathds 1}}

\def\e{{\varepsilon}}

\def\ba{{\backslash}}

\def\D{\Delta}

\def\sm{\setminus}
\def\ba{\setminus}

\def\d{\delta}

\def\e{\varepsilon}

\def\phi{\varphi}

\def\g{\gamma}

\def\k{\kappa}

\def\r{\rho}

\def\s{\sigma}

\def\o{\omega}

\def\D{\Delta}

\def\L{\Lambda}

\def\G{\Gamma}

\def\O{{\Omega}}

\def\T{\T}

\def\SS{{\mathcal{F}}}

\def\GG{{\mathcal{G}}}
\def\PP{{\mathcal P}}

\def\V|{{\Vert}}


\begin{document}

\author{Benedikt Jahnel}
\address[Benedikt Jahnel]{Weierstrass Institute Berlin, Mohrenstr. 39, 10117 Berlin, Germany, \texttt{https://www.wias-berlin.de/people/jahnel/}}
\email{Benedikt.Jahnel@wias-berlin.de}

\author{Christof K\"ulske}
\address[Christof K\"ulske]{Ruhr-Universit\"at   Bochum, Fakult\"at f\"ur Mathematik, D44801 Bochum, Germany, \texttt{http://www.ruhr-uni-bochum.de/ffm/Lehrstuehle/Kuelske/kuelske.html}}
\email{Christof.Kuelske@ruhr-uni-bochum.de}


\title{Attractor properties for irreversible and reversible interacting particle systems}

\date{\today}

\maketitle

\begin{abstract}
We consider translation-invariant interacting particle systems on the lattice with finite local state space admitting at least one Gibbs measure as a time-stationary measure. The dynamics can be irreversible but should satisfy some mild non-degeneracy conditions. We prove that weak limit points of any trajectory of translation-invariant measures, satisfying a non-nullness condition, are Gibbs states for the same specification as the time-stationary measure. This is done under the additional assumption that zero entropy loss of the limiting measure w.r.t.~the time-stationary measure implies that they are Gibbs measures for the same specification. 

We show how to prove the non-nullness for a large number of 
cases, and
also give an alternate version of the last condition such that the non-nullness requirement can be dropped. 
As an application we obtain the attractor property if there is a reversible Gibbs measure. 

Our method generalizes convergence results using relative entropy techniques to a large class of dynamics including irreversible and non-ergodic ones.
\end{abstract}

\section{Introduction}
The last years have seen an interest in the analysis of infinite-volume lattice measures 
under stochastic time-evolutions, with a particular view on the possible production of singularities of such measures  \cite{EnFeHoRe02,EnFeHoRe10,EnRu09,ErKu10,KuNy07,KuRe06,KuLeRe04}. 
For analogues of this phenomenon for spatially structured systems 
beyond the lattice world see also~\cite{JaKu16}.
These singularities are related to the emergence of long spatial memory in the conditional probabilities of the time-evolved 
measures at given transition times. When the initial measure is a Gibbs measure in a low-temperature phase 
for some absolutely summable potential it may happen that a time-evolved potential ceases to exist, and one speaks of a Gibbs-non-Gibbs transition. 
Such phenomena are proved to occur on the lattice for weakly interacting Glauber dynamics, based on the detection 
of 'hidden phase transitions'.  Suggested by mean-field analogues, singularities are 
expected to appear (and even more easily so) for strongly interacting reversible dynamics.  
While the focus of this research has been much on reversible dynamics, one expects
similar singularities during time-evolution in the huge field of irreversible dynamics, see for example \cite{Ma06,Li85}, which are even harder to analyze.  

While the assumption that all time-evolved measures 
are Gibbs would make it easy to obtain a good control in terms of finite-volume 
approximations, 
this possible occurrence of non-localities in turn poses difficulties to control the large-time behavior of trajectories of time-evolved measures \cite{JaKu15}.  
It is the purpose of this paper to exploit the key concept of relative entropy change 
(per site)
along 
trajectories in this context, including situations of multiple phases, and including situations of irreversible dynamics. As it turns out this lack of reversibility 
forces us to work much harder, and is responsible for opening 
our discussion of the non-nullness property. 

Relative entropy has a huge importance in the probability theory of statistical mechanics in infinite volume, via its relevance in large deviations, via the 
Gibbs variational principle, see for example \cite{Ge11, KuLeRe04, De15}, 
and also via a new formulation to analyze Gibbs-non Gibbs transitions in terms of a variational principle in path space, see\cite{EnFeHoRe10, ErKu10}. 
The study of entropy decay to unique equilibria connects probability, analysis, and geometry in fascinating ways, see \cite{Vi09,ReSt05,ErKuSt15,BaGeLe14}.
Its successful use as a Lyapunov function in infinite volume in the context of stochastic time-evolutions goes back to very early work of Holley \cite{Ho71,HoSt77} in the special case 
of a Glauber dynamics 
for Ising models. 
Related methods are used in  the excellent reference 
\cite{Ge79} for dynamics with conserved particle numbers, 
however under the necessary assumption of reversibility. 
Already in \cite{Ku84}, zero entropy loss is used to classify invariant states, but the more 
difficult issue of the behavior of trajectories 
for starting measures off the invariant states is not studied. 
In this paper we build up on these initial steps and go  beyond the reversible case. 
As an application we 
also provide a treatment of general Glauber dynamics in infinite volume, with finite alphabets
and general Hamiltonians 
(which to our knowledge has not appeared 
in the literature).   

\medskip
We work in the setting of stochastic dynamics for lattice systems in the infinite volume and in continuous time. 
Our local state spaces are finite and the dynamics is specified by giving the rates to jump between 
different symbols in this alphabet. 
In most interesting cases these jump processes are non-independent over the sites which creates the possibility for macroscopically non-trivial collective behavior. 
In all what follows we assume lattice-translation invariance for the rules specifying the 
dynamics.
We will for the most part {\em not assume} that the dynamics is reversible for a particular measure. 

We look at initial configurations which are chosen from lattice-translation invariant starting measures, 
and will then be interested in the corresponding trajectories of lattice-translation invariant infinite-volume measures. 
We ask for possible large time limiting behavior. In the language of dynamical systems, we want to know 
the omega-limit sets of the dynamics, that is the set of possible weak limit points of $\nu_{t_n}$ where $t_n$ tends to infinity. Here the usual weak convergence is chosen in which 
convergence of measures is checked in terms of local observables. In particular, by compactness, 
there are always weak limit points.  
The dynamics has at least one time-stationary measure $\mu$, which might 
be ergodic w.r.t.~lattice translations or not. In fact we have examples for both situations.  
For this measure we will assume that it is even a Gibbs measure w.r.t.~a quasilocal specification $\g$, in other words $\mu\in\GG(\g)$.  This is the case for Glauber dynamics, and for a  class of irreversible dynamics \cite{JaKu12}. However, there are examples of irreversible dynamics with non-Gibbsian invariant measures \cite{EnFeSo93}. If there is one non-Gibbsian invariant measure, the other invariant measures must be non-Gibbsian too \cite{Ku84}. 
We give an example of non-Gibbsian invariant measures caused by the lack 
of reducibility of the dynamics in Subsection \ref{fuenfdrei}.

Let us mention that  relative entropy techniques have 
been used successfully also in related but different ways:
In the hydrodynamic-limit approach \cite{KiLa99} convergence of an interacting particle system to some PDE system under the thermodynamic limit is studied for finite times.
In this context, the so-called relative entropy method is based on the observation (see \cite{Ya91,OlVa91,OlVaYa91}) that, if the relative entropy density of some starting measure is zero w.r.t.~the limiting measure driven by the PDE, then this is also true for the time-evolved measure. In that sense the technique does not directly tackle the question of large-time behaviour of the original system. In contrast to these hydrodynamic limit statements, see for example \cite{Sp93, FuSp97}, which are based on relative entropy techniques, in our setting, there is no spatial rescaling involved.

\medskip
Here we want to use the relative entropy $h(\nu|\mu)$ w.r.t.~to a time-stationary measure $\mu\in \GG(\g)$
as a Lyapunov function to investigate trajectories and 
limit points. More precisely, we give criteria under which the set of weak limit points is contained in the set of Gibbs measures associated to the invariant measure.

Let us note that in the particular case of the uniqueness regime, $|\GG(\g)|=1$, the subject of entropy decay under time-evolutions is intimately linked to Log-Sobolev inequalities for infinite-volume measures, see \cite[Chapter 5]{BaGeLe14} or \cite{GuZe03}. 
Proving a Log-Sobolev inequality for a non-equilibrium model implies the exponential decay of the relative entropy distance and thus gives not only the attractor property but also the rate of convergence 
to the unique equilibrium. 
We cannot use these methods here since our interest goes beyond situations of uniqueness to situations where multiple invariant measures may occur. 
Instead, our method is based on semicontinuity of the entropy {\em loss} functional 
in infinite volume. 

\subsection{Relative entropy as a Lyapunov function}
The difficulties using the relative entropy as a Lyapunov function in the infinite volume are caused by the potential lack of continuity. 
Recall that the  
relative entropy density $\nu\mapsto h(\nu|\mu)$ 
is a lower semicontinuous function in the weak topology, 
but in general it is not upper semicontinuous. 
Looking at the time-derivative of the relative entropy $g_L(\nu_{t_n}|\mu)$, as defined in \eqref{REL}, along trajectories which 
are sampled at time instances $t_n$ tending to infinity,  
we have $\lim_{n\uparrow\infty}g_L(\nu_{t_n}|\mu)= 0$. 
We would like to conclude that $g_L(\nu^*|\mu)= 0$ where $\nu^*$ denotes a weak limit point of the trajectory.
This equation expresses zero entropy loss of the limiting measure and 
is in itself very useful to characterize possible limits. In interesting cases it may have 
multiple solutions $\nu^*$.
In many cases, and even irreversible situations as in \cite{JaKu12}, it allows to characterize $\nu^*$ by concluding 
that these solutions must be elements of $\GG(\g)$.

Now, in order to prove $g_L(\nu^*|\mu)= 0$, 
we would have to know that $\nu \mapsto g_L(\nu|\mu)$ is upper semicontinuous.
A proof that  $g_L$ is upper semicontinuous has been given in a reversible situation for the particular case of 
the stochastic Ising model  for which the corresponding 
Ising Gibbs measures are reversible measures in \cite{Ho71,HiSh75}. It is the prime aim of the present paper to move into the realm of non-reversibility.
As our main result, we prove that $g_L$ is upper semicontinuous also for general types of 
non-reversible dynamics 
and also in the situation of general finite state spaces. 
As a byproduct we also prove the attractor property for reversible dynamics w.r.t.~Gibbs measures for irreducible finite state space interacting particle systems (IPS) on the lattice.

\subsection{A useful decomposition of entropy loss} 
We are looking for monotonicity in suitably chosen finite-volume approximations
of $g_L$ to conclude that $g_L$ is upper semicontinuous. 
This is what is done in \cite{HiSh75}, making heavy use of {\em reversibility} which we cannot use in general,  
and taking advantage 
of the two-spin situation. 
We work with a useful decomposition of $g_L$ to treat also irreversible dynamics, 
see \eqref{RelEntLoss} and \eqref{EntropyChange1zt}.  This decomposition is explained for 
single-site dynamics in Appendix~\ref{AP}.
The decomposition does not have an interpretation in terms 
of decomposition of generators into symmetric and antisymmetric part. 
It is made to separate dangerous terms 
in such a way as to bring volume monotonicity and  convexity into play to help provide semicontinuity. 
That it can be made to work in infinite volume as well, is a main result of the present paper, and relies 
on the control of boundary terms. Here irreversibility poses additional difficulties compared to for example Glauber dynamics, and the models considered in~\cite{Ge79}, and the property of non-nullness along the trajectory becomes 
essential. 
We also give a proof of this non-nullness property along the trajectory 
for finite-range dynamics with single-site updates, which may be irreversible.
We stress that we do not need any assumption on quasilocality along the trajectory  
which in many cases indeed would not hold, as non-Gibbsian measures are known 
to occur under stochastic dynamics, even under independent dynamics, see for example \cite{EnFeHoRe02,EnFeHoRe10,EnRu09,ErKu10,KuNy07,KuRe06,KuLeRe04}.

\medskip
One particular motivation for considering the relative entropy decay under irreversible dynamics comes from a class of models we consider in \cite{JaKu12,JaKu14b}. These models exhibit dynamical non-ergodicity, in the sense of IPS, in the presence of a unique time-stationary Gibbs measure, making rigorous a heuristics of \cite{MaSh11}. In the analysis of a mean-field version of these rotation dynamics in \cite{JaKu14} we were able to show the attractor property of the limiting cycle using relative entropy techniques on finite-dimensional simplexes. This proves synchronization in the sense of attractivity of macroscopically coherent rotating states. Let us mention that a similar type of synchronization is also frequently studied for other, however mostly mean-field models, for example the Kuramoto model for coupled noisy phase oscillators \cite{GPP12}.
Beyond that, our results provide a very general analytical approach applicable also to systems where coupling and duality 
tools are not available.

\subsection{Organisation of the manuscript}
In Section~\ref{Equ} we first introduce the equilibrium setting of infinite-volume Gibbs measures and relative entropy densities. In Section~\ref{IPSSetting} we present the dynamical setting of interacting particle systems and the associated relative entropy loss densities. In Sections~\ref{AttrIrr}, \ref{AttrIrr2} and~\ref{AttrRe} the main result about the attractor property for irreversible and reversible dynamics are stated. In Section~\ref{Sup} we present supporting results which allow us to prove the main theorems of the paper. The technical proofs of the supporting results are contained in Section~\ref{Proofs}. Finally, in the Appendix~\ref{AP} we discuss
the special case of independent dynamics.

\subsection{Acknowledgments}
The authors thank the editor and anonymous referees for comments and suggestions that helped to improve the presentation of the material. This research was supported by the Leibniz program Probabilistic Methods for Mobile Ad-Hoc Networks.

\section{Entropy decay for interacting particle systems}\label{Sec2}

\subsection{Gibbs measures and relative entropy}\label{Equ}
Let $\PP_\theta$ denote the set of translation-invariant probability measures on the configuration space $\O=\{1,\dots,q\}^{\Z^d}$ equipped with the usual product topology and the corresponding Borel sigma-algebra $\SS$. Then, for $\mu, \nu\in\PP_\theta$ and a finite set of sites $\L\Subset\Z^d$ define the \textit{relative entropy} via
\begin{equation*}\label{LRE}
\begin{split}
h_\L(\nu|\mu):=
\begin{cases}
  \sum_{\o_\L\in\{1,\dots,q\}^\L}\nu(\one_{\o_\L})\log\frac{\nu(\one_{\o_\L})}{\mu(\one_{\o_\L})}
,  & \text{if } \nu\ll\mu,\\
  \infty, & \text{else, }
\end{cases}
\end{split}
\end{equation*}
where $\one_{\o_\L}$ denotes the indicator function on the finite configuration $\o_\L$, and $\nu(f)=\int\nu(d\o)f(\o)$ is a short-hand notation for integration. Further, define
the \textit{relative entropy density} via
\begin{equation*}\label{SRE}
\begin{split}
h(\nu|\mu):=\lim_{n\uparrow\infty}\frac{1}{|\L_n|}h_{\L_n}(\nu|\mu)
\end{split}
\end{equation*}
where $\L_n:=[-2^n+1,2^n-1]^d$ is a sequence of hypercubes centered at the origin, whenever the limit exists. 

\medskip
We will be interested in situations where $\mu$ is a Gibbs measure for a translation-invariant non-null quasilocal specification on $\O$. A \textit{specification} is a family $\g =(\g_\L)_{\L\Subset\Z^d}$ of proper probability kernels $\g_\L(\one_{\eta_\L}|\eta_{\L^c})$ from $(\O,\SS_{\L^c})$ to the set of probability measures on $(\O,\SS)$. Here, $\SS_{\L^c}$ is the sub-sigma algebra of $\SS$ generated by the open sets in $\{1,\dots,q\}^{\L^c}$ and properness means that if $\D\subset\L^c$, then $\g_\L(\one_{\eta_\L}\one_{\eta_\D}|\eta_{\L^c})=\g_\L(\one_{\eta_\L}|\eta_{\L^c})\one_{\eta_\D}(\eta_{\L^c})$. Further, specifications satisfy the consistency condition $\g_\L(\g_\D(\one_{\eta_\D}|\cdot)|\eta_{\L^c})=\g_\L(\one_{\eta_\D}|\eta_{\L^c})$ whenever $\D\subset\L$. 

\medskip
In the following, we will often denote for a given configuration $\o\in\O$ by $\o_\L$ its projection to the volume $\L\subset\Z^d$ and to write $\o_\L\o_\D$ for the finite-volume configuration in $\L\cup\D$ composed of $\o_\L$ and $\o_\D$ with disjoint $\L,\D\subset\Z^d$. We also denote $\L^c:=\Z^d\sm \L$ and write $i^c$ instead of $\{i\}^c$ for $i\in \Z^d$.
\begin{defn}\label{Spec}
The specification $\g$ is called 
\begin{enumerate}
\item \emph{translation invariant}, if for all $\L\Subset\Z^d$ and $i\in\Z^d$ we have $\g_{\L+i}(\one_{\eta_{\L+i}}|\eta_{{(\L+i)}^c})=\g_\L(\one_{\eta_\L}|\eta_{\L^c})$ where $\L+i$ denotes the lattice translate of $\L$ by $i$;
\item \emph{non-null}, if $\inf_{\eta}\g_o(\one_{\eta_o}|\eta_{o^c})\ge\d$ for some $\d>0$;
\item \emph{quasilocal}, if for all $\L\Subset\Z^d$, $\lim_{\D\uparrow\Z^d}\sup_{\eta,\xi}|\g_\L(\one_{\eta_\L}|\eta_{\D\setminus\L}\xi_{\D^c})-\g_\L(\one_{\eta_\L}|\eta_{\L^c})|=0$.
\end{enumerate}
\end{defn}
The infinite-volume probability measure $\mu$ is called a \textit{Gibbs measure} for $\g$, i.e., $\mu\in\GG(\g)$, if $\mu$ satisfies the DLR equation, namely for all $\L\Subset\Z^d$ and $\eta_\L$ we have $\mu(\g_\L(\one_{\eta_\L}|\cdot))=\mu(\one_{\eta_\L})$. For details on Gibbs measures and specifications see \cite{Ge11,EnFeSo93}.

\medskip
In order to guarantee existence of the relative entropy density, $\mu$ has  to be \textit{asymptotically decoupled} as defined in \cite{Pf02,KuLeRe04}. For this, denote $\L_n$ the centered box with side-length $2n+1$. 
\begin{defn}\label{AsDeMe}
A probability measure $\mu$ on $(\O,\SS)$ is called \emph{asymptotically decoupled} if
\begin{enumerate}
\item there exist $d:\N\mapsto\N$ and $c:\N\mapsto[0,\infty)$, such that 
$$\lim_{n\uparrow\infty}d(n)/n=0 \text{ and }\lim_{n\uparrow\infty}c(n)/|\L_n|=0.$$
\item for all $i\in\Z^d$, $n\in\N$, $A\in\SS$ measurable w.r.t.~$\L_n+i$ and $B\in\SS$ measurable w.r.t.~$(\L_{n+d(n)}+i)^c$, we have 
$$e^{-c(n)}\mu(A)\mu(B)\le\mu(A\cap B)\le e^{c(n)}\mu(A)\mu(B).$$
\end{enumerate}
\end{defn}

The following result, proved in \cite[Proposition 3.2]{Pf02}, guarantees existence of the relative entropy density w.r.t.~asymptotically decoupled measures.
\begin{lem}\label{RED_Exists}
Let $\nu,\mu\in\PP_\theta$ and $\mu$ asymptotically decoupled. Then, the relative entropy density $h(\nu|\mu)$ exists and is non-negative.
\end{lem}
For example specifications defined via translation-invariant uniformly 
absolutely summable potentials $\Phi=(\Phi_A)_{A\Subset\Z^d}$ are translation invariant, non-null and quasilocal. Gibbs measures for such \textit{Gibbsian specifications} are moreover asymptotically decoupled and hence the relative entropy density of any translation-invariant measure
relative to them exists. 

\medskip
Note, as a subtlety, that for general translation-invariant specifications without any further assumptions on locality properties, existence of an absolutely summable \textit{translation-invariant} potential is not guaranteed, see \cite{Su73,Ko74,EnFeSo93}. This is why we are imposing asymptotic decoupledness as an additional requirement. As a side-remark, for systems of point particles going back and forth between 
specifications and potential representations are even more subtle, 
see~\cite{JaKu17}.
The equilibrium model considered in \cite{JaKu12} provides an example of such asymptotically decoupled Gibbs measures, where the specification is a priori not given in terms of an absolutely summable translation-invariant potential. 

\subsection{IPS dynamics and the relative entropy loss}\label{IPSSetting}
Consider time-continuous, translation-invariant Markovian dynamics on 
$\O$, namely IPS characterized by time-homogeneous generators $L$ with domain $D(L)$ and its associated Markovian semigroup $(P^L_t)_{t\geq0}$. For the IPS we adopt the exposition given in 
the standard reference
\cite[Chapter I]{Li85}. In all generality the generator $L$ is given via jump-rates $c_\D(\eta,\xi_\D)$ in finite volumes $\D\Subset\Z^d$, continuous in the starting configurations $\eta\in\O$
\begin{equation}\label{IPS}
Lf(\eta)=\sum_{\D\Subset\Z^d}\sum_{\xi_\D}c_\D(\eta,\xi_\D)[f(\xi_\D\eta_{\D^c})-f(\eta)].
\end{equation}
To ensure well-definedness, the jump-measures must satisfy a number of conditions. Most importantly the single-site jump-intensities have to be bounded, i.e., for
$c_\D(\eta):=\sum_{\xi_\D\neq\eta_\D}c_\D(\eta,\xi_\D)$ and $c_\D:=\sup_{\eta}c_\D(\eta)$
we assume $\sum_{\D\ni o}c_\D<\infty$. In fact the definition of $L$ in \eqref{IPS} should be read in such a way that the two summations are only over those $\D$ and $\xi_\D$ with $c_\D(\eta,\xi_\D)>0$. We will call an IPS well-defined if it is well-defined in the sense of \cite[Chapter I]{Li85}.

The following additional conditions on IPS will be used in the sequel. 
\begin{defn}\label{Gibbs-Attractor-IPS}
Let $Lf(\eta)=\sum_{\D\Subset\Z^d}\sum_{\xi_\D}c_\D(\eta,\xi_\D)[f(\xi_\D\eta_{\D^c})-f(\eta)]$ 
be a well-defined IPS. We say that 
\begin{enumerate}
\item\label{tran} $L$ is \emph{translation invariant}, if all rates are translation invariant;
\item\label{fini} for $L$ there are only \emph{finitely many types of transitions}, if $c_\D>0$ for only finitely many $\D\Subset\Z^d$;
\item\label{cont} for $L$ the rates are \emph{uniformly continuous}, if 
$\lim_{\L\uparrow\Z^d}\sup_{\D\ni o}\sup_{\eta,\xi,\s}|c_\D(\eta_{\L}\s_{\L^c},\xi_\D)-c_\D(\eta,\xi_\D)|=0$;
\item\label{mini} $L$ has a \emph{strictly positive minimal transition rate}, if for all $\D\Subset\Z^d$ with $c_\D>0$
we have 
$\inf_{\eta,\xi:c_\D(\eta,\xi_\D)>0}c_\D(\eta,\xi_\D)>0$;
\item\label{trap} $L$ can \emph{not enter trap states}, if for all $\eta\in\{1,\dots,q\}^{\Z^d}$ and $\xi_\D\in\{1,\dots,q\}^{\D}$ we have that $c_\D(\eta,\xi_\D)>0$ implies that $c_\D(\xi_\D\eta_{\D^c})>0$.
\end{enumerate}
\end{defn}

Examples of IPS satisfying the above conditions are the stochastic Ising model or more general Glauber dynamics. As another example consider the exclusion process on $\{0,1\}^{\Z^d}$ with rates 
$$c_{\{x,y\}}(\eta;(1-\eta_x),(1-\eta_y))=p(x,y)\eta_x(1-\eta_y)+p(y,x)\eta_y(1-\eta_x).$$
Here, $p(x,y)$ describes the possibly non-symmetric rate of moving a particle from $x$ to site $y$. Such processes are contained in the class of IPS satisfying above conditions as long as $p(x,y)>0$ implies $p(y,x)>0$. Contact processes or voter models have trap states and thus our approach can not be applied.

\medskip
In this paper we want to analyze Gibbsian models given in terms of translation-invariant non-null quasilocal specifications $\g$ and transformations given by translation-invariant IPS dynamics which have at least one of the Gibbs measures as a time-stationary measure. Our main tool is to consider the evolution of the relative entropy density. Let us define for any $\nu,\mu\in\PP_\theta$ with $\mu P_t^L=\mu$ and $\L_n$, the \textit{relative entropy loss} via
\begin{equation*}\label{RelEntLos}
\begin{split}
g_L^n(\nu|\mu):=\frac{d}{dt}_{|t=0}h_{\L_n}(\nu P_t^L|\mu).
\end{split}
\end{equation*}
Similar as for the relative entropy, we define the \textit{relative entropy loss density} via
\begin{equation}\label{REL}
\begin{split}
g_L(\nu|\mu):=\limsup_{n\uparrow\infty}\frac{1}{|\L_n|}g_L^{\L_n}(\nu|\mu)
\end{split}
\end{equation}
where $\L_n=[-2^n+1,2^n-1]^d$. In Proposition~\ref{EntropyLossLem} we show that $g_L(\nu|\mu)\le0$, which justifies the name of $g_L(\nu|\mu)$.

\medskip
In order for our main results to cover also models with conserved quantities, let us introduce the notation $D_L(\mu)$ for a weakly-closed subset of $\PP_\theta$ which is time stationary under the dynamics $L$ containing a particular measure $\mu\in\PP_\theta$, see also below for some specific examples. In the following section we state our first main result about the attractor property for not necessarily reversible IPS.

\subsection{Attractor property for irreversible dynamics with non-nullness}\label{AttrIrr}
For our first main result, we will assume that under the dynamics the following \emph{zero entropy loss condition} holds. 
\begin{cond}\label{zero entropy loss condition}
We say that a well-defined IPS dynamics $L$ satisfies the \emph{zero entropy loss condition} if the following is true.
\begin{enumerate}
\item $L$ satisfies the conditions~(\ref{tran}), (\ref{fini}), (\ref{cont}), (\ref{mini}) and (\ref{trap}) in Definition~\ref{Gibbs-Attractor-IPS}. 
\item For $L$ there exists a translation-invariant asymptotically-decoupled time-stationary $\mu\in\GG(\g)$ where $\g$ is translation-invariant, non-null and quasilocal.
\item For any $\nu\in D_L(\mu)$ with $g_L(\nu|\mu)=0$ it follows that $\nu\in\GG(\g)$.
\end{enumerate}
\end{cond}

Without the time derivative, i.e., with $g_L(\nu|\mu)$ replaced by $h(\nu|\mu)$, and for $D_L(\mu)=\PP_\theta$, this condition is one direction of the Gibbs variational principle, see for example \cite[Theorem 15.37]{Ge11}. 
For $D_L(\mu)=\PP_\theta$, all conditions given in Definition \ref{Gibbs-Attractor-IPS} plus the above Condition \ref{zero entropy loss condition} involving the time-derivative are proved to hold for example for the stochastic Ising model in~\cite{Li85,Ho71,HoSt77} for any choice of a possibly non-unique translation-invariant equilibrium measures $\mu$. It can also be proved to hold for non-reversible rotation dynamics see~\cite{JaKu12}. 

A crucial requirement for our first main result to hold is \emph{non-nullness}.
\begin{defn}\label{Def_Non_Null}
A probability measure $\nu$ on $\O$ is \emph{non-null} if there exists $\d>0$ and a version of the single-site conditional probabilities such that $\nu(\one_{\eta_o}|\eta_{o^c})\ge\d$ for $\nu$-a.a.~$\eta$.
\end{defn}
 Time-evolved random fields should be non-null under rather weak assumptions on the dynamics. We include a full proof of the following statement on non-nullness for finite-range single-site generators.
\begin{prop}\label{Prop_Non_Null}
Consider a translation-invariant IPS generator with single-site updates of the form 
$$Lf(\eta)=\sum_{i\in \Z^d}\sum_{\xi_i=1,\dots,q}c_i(\eta,\xi_i)[f(\xi_i\eta_{i^c})-f(\eta)].$$
We further assume that 
\begin{enumerate}
\item $L$ is finite range. That is, there exists a finite centered box $\D\Subset\Z^d$ such that for all $\xi_o$, we have $c_o(\eta,\xi_o)=c_o(\eta_{\D}\s_{\D^c},\xi_o)$ for all $\eta,\s\in\O$. 
\item for $L$, reachability is independent of the boundary conditions. That is, 
whenever $c_o(\eta,\xi_o)>0$ also $c_o(\s,\xi_o)>0$ for all $\s$ 
with $\s_o=\eta_o$. 
In this case we say that $\xi_o$ can be reached from $\eta_o$, and write $d(\eta_o,\xi_o)$ for the indicator of this event. 
\item $L$ is single-site irreducible. That is, the Markov chain on the single state space 
with rates given by $d$ is irreducible.
\end{enumerate}
Then, for all $\tau>0$ there exists a $\d=\d(\tau)>0$ such that for any starting measure $\nu\in\PP_\theta$, and any time $t\in[\tau,\infty)$, the 
time-evolved measure $\nu P^L_t$ is non-null with constant $\d$. Furthermore, any subsequential weak limiting measure $\nu_*=\lim_{n\uparrow\infty}\nu P_{t_n}^L$ where $t_n\uparrow\infty$ is also non-null with constant $\d$.
\end{prop}
Let us note that the above result allows for irreversible dynamics. As it will become clear from 
the proof, the statement can be generalized also to multi-site updates and rates that are not strictly finite range. 
We can now state our first main result. 
\begin{thm}\label{MainCor}
Let the well-defined IPS dynamics $L$ satisfy Condition~\ref{zero entropy loss condition} with time-stationary $\mu\in\GG(\g)$. Then, for any translation-invariant starting measure $\nu\in D_L(\mu)$ where the sequence $(\nu P_{t_n}^L)_{n\in\N}$ consists of non-null probability measures and converges weakly to the non-null probability measure $\nu_*$ as $t_n\uparrow\infty$, we have that $\nu_*\in\GG(\g)$.
\end{thm}

In the next section we show that the non-nullness requirement can be dropped if the zero entropy loss Condition \ref{zero entropy loss condition} is replaced by an approximating zero entropy loss condition. 
\subsection{Attractor property without non-nullness assumption}\label{AttrIrr2}
Consider another sequence of centered hypercubes 
$\tilde\L_n=[-2^n+n,2^n-n]^d\subset\L_n= [-2^n+1,2^n-1]^d$. We then define the \textit{approximating relative entropy loss}
as
\begin{align*}
\tilde g_L^n(\nu|\mu)
=\sum_{i\in\tilde\L_n}\sum_{\D\ni i}\tfrac{1}{|\D\cap\L_n|}\sum_{\xi_{\D}}\int\nu(d\eta)c_\D(\eta_{\L_n},\xi_{\D})\log\tfrac{\nu(\one_{\xi_{\D\cap\L_n}}|{\eta_{\L_n\setminus\D}})\mu(\one_{\eta_{\D\cap\L_n}}|{\eta_{\L_n\setminus\D}})}{\nu(\one_{\eta_{\D\cap\L_n}}|{\eta_{\L_n\setminus\D}})\mu(\one_{\xi_{\D\cap\L_n}}|{\eta_{\L_n\setminus\D}})}.
\end{align*}
This definition seems technical at first reading, but in view of the representation of $g_L(\nu|\mu)$ given in Lemma~\ref{Lem_Rep} it becomes clear, that $\tilde g_L(\nu|\mu)$ is equivalent to $g_L(\nu|\mu)$ except that updates are only performed on those sets $\D$ that communicate with the smaller volume $\tilde\L_n$. As before, we define the \textit{approximating relative entropy loss density} via
\begin{equation*}\label{Gee}
\begin{split}
\tilde g_L(\nu|\mu):=\limsup_{n\uparrow\infty}\frac{1}{|\L_n|}\tilde g_L^{n}(\nu|\mu).
\end{split}
\end{equation*}
Let us assume that under the dynamics the following \textit{approximating zero entropy loss condition} holds:
\begin{cond}\label{approx zero entropy loss condition}
We say that a well-defined IPS dynamics $L$ satisfies the \emph{approximating zero entropy loss condition} if the following is true.
\begin{enumerate}
\item $L$ satisfies the conditions~(\ref{tran}), (\ref{fini}), (\ref{cont}), (\ref{mini}) and (\ref{trap}) in Definition~\ref{Gibbs-Attractor-IPS}. 
\item For $L$ there exists a translation-invariant asymptotically-decoupled time-stationary $\mu\in\GG(\g)$ where $\g$ is translation-invariant, non-null and quasilocal.
\item For any $\nu\in D_L(\mu)$ we have $\tilde g_L(\nu|\mu)\le 0$ and from $\tilde g_L(\nu|\mu)=0$ follows that $\nu\in\GG(\g)$.
\end{enumerate}
\end{cond}
In case of the SEP, also the Condition~\ref{approx zero entropy loss condition} with $D_L(\mu)$ given by the set of translation-invariant probability measures with particle density as $\mu$ can be verified. Note that $D_L(\mu)$ is indeed weakly closed as it is defined in terms of the expectation of the spin variable at the origin. 
Another example for which the Condition~\ref{approx zero entropy loss condition} can be verified with $D_L(\mu)=\PP_\theta$ is the stochastic Ising model, see~\cite{HiSh75}.
Under the approximating zero entropy loss condition we can prove the attractor property avoiding a non-nullness requirement for the trajectory and the limiting measure.

\begin{thm}\label{MainCorApprox}
Let the well-defined IPS dynamics $L$ satisfy Condition~\ref{approx zero entropy loss condition} with time-stationary $\mu\in\GG(\g)$. Then, for any translation-invariant starting measure $\nu\in D_L(\mu)$ where the sequence $(\nu P_{t_n}^L)_{n\in\N}$ converges weakly to $\nu_*$ as $t_n\uparrow\infty$, we have that $\nu_*\in\GG(\g)$.
\end{thm}
In case of the SEP, Theorem~\ref{MainCorApprox} provides an alternative proof of the well known attractor result, see~\cite{Li85}, via relative entropy.
Why do we need to pay attention to the difference between $\tilde g_L(\nu|\mu)$ and $g_L(\nu|\mu)$? We introduce the approximating quantity $\tilde g^n_L(\nu|\mu)$ to offer a way to circumvent the non-nullness requirement. At the heart of most of the arguments in this paper is the need to bound certain boundary terms in $\L_n\setminus\tilde\L_n$. This is difficult to do in general. We offer two solutions which do the job from two different perspectives: 
\begin{enumerate}
\item Ensure that the sequence of measures $\nu_{t_k}$ is non-null. This can be done sometimes, see Proposition 2.7.
\item Use a different Lyapunov function, which is $\tilde g$ instead of $g$. Then no non-nullness is needed, but one has to verify that it really can serve as a Lyapunov function, which is what we require in Condition~\ref{approx zero entropy loss condition}. This is also the strategy used in our reference results, for example~\cite{HiSh75}, for the stochastic Ising model. This is also the strategy we use to prove Theorem~\ref{GTilde} presented in the next section.
\end{enumerate}

In the following section, we prove that the approximating zero entropy loss condition is satisfied if the time-stationary measure $\mu$ is even \emph{reversible} for $L$. Hence, as an application of Theorem~\ref{MainCorApprox}, the attractor property is proven for reversible dynamics.

\subsection{Attractor property for reversible dynamics}\label{AttrRe}
For our final result, it suffices to show that Condition~\ref{approx zero entropy loss condition}
can be verified if $\mu$ is a reversible measure for $L$ and the requirement that $L$ has no trap states is replaced by the following stronger assumption of irreducibility.
\begin{defn}\label{Gibbs-Attractor-IPS_Irre}
Let $Lf(\eta)=\sum_{\D\Subset\Z^d}\sum_{\xi_\D}c_\D(\eta,\xi_\D)[f(\xi_\D\eta_{\D^c})-f(\eta)]$ 
be a well-defined translation-invariant IPS. We say that $L$
is irreducible, if for all $\eta^{(0)}\in\O$ and $\s\in\{1,\dots,q\}^{\D}$ with $\D\Subset\Z^d$ there exists a finite sequence of configurations $\{\eta^{(1)},\dots,\eta^{(n)}\}$ with $\eta^{(i)}\in\O$ and  $\eta^{(n)}=\eta_{\D^c}^{(0)}\s_\D$ such that the transition rates to jump from $\eta^{(i-1)}$ to $\eta^{(i)}$ are positive for all $i\in\{1,\dots,n\}$.
\end{defn}
Note that for example Kawasaki dynamics or the SEP are not irreducible in the above sense.
The following Theorem~\ref{GTilde} is our main statement about 
the attractor property for reversible dynamics and generalizes~\cite{Ho71, HiSh75}.

\begin{thm}\label{GTilde}
Let the well-defined IPS dynamics $L$ satisfy the conditions~(\ref{tran}), (\ref{fini}), (\ref{cont}) and (\ref{mini}) in Definition~\ref{Gibbs-Attractor-IPS}. 
Further, assume that $L$ is irreducible and admits 
a reversible Gibbs measure $\mu\in\GG(\g)$ which is translation-invariant, asymptotically-decoupled and where $\g$ is translation-invariant, non-null and quasilocal. 
Then, for any $\nu\in\mathcal{P}_\theta$ where the sequence $(\nu P_{t_n}^L)_{n\in\N}$ converges weakly to $\nu_*$ as $t_n\uparrow\infty$, we have that $\nu_*\in\GG(\g)$.
\end{thm}
In the next section we present the strategy of proof, state the supporting results and prove the main results for the irreversible dynamics. The proofs of the supporting results follow in Section~\ref{Proofs}. In the Appendix~\ref{AP} we also present a comprehensive single-site example to illustrate our general strategy. 

\section{Strategy of proof}\label{Sup}
The strategy of the proof, in a nutshell, is given by the following sequence of inequalities
\begin{equation}\label{nutshell}
\begin{split}
0&=\lim_{k\uparrow\infty}g_L(\nu_{t_k}|\mu)\le g_L(\nu_*|\mu)\le 0
\end{split}
\end{equation}
which, under the zero entropy loss Condition, then implies Gibbsianness of the limiting measure. All of the above inequalities require a proper introduction and this is what we are going to do now. Let us start with the last inequality in \eqref{nutshell} and state that indeed the relative entropy is non increasing under the time evolution. 
\begin{prop}\label{EntropyLossLem}
Let $\nu\in\PP_\theta$ and $L$ be a well-defined IPS generator satisfying conditions~(\ref{tran}), (\ref{fini}) and (\ref{cont}) in Definition~\ref{Gibbs-Attractor-IPS}.
Assume that for $L$ there exists a translation-invariant asymptotically-decoupled time-stationary Gibbs measure $\mu\in\GG(\g)$ where $\g$ is translation-invariant non-null and quasilocal, then $g_L(\nu|\mu)\le \tilde g_L(\nu|\mu)$. If $\nu$ is additionally non-null, then $g_L(\nu|\mu)\le0$.
\end{prop}
The relative entropy density w.r.t.~translation-invariant probability measures is a non-increasing function under rather general transformations, see for example \cite[Lemma 3.3]{EnFeSo93}. In case of dynamics with conserved particle numbers and allowing for a Gibbs measure as a \emph{reversible} measure, this statement is proved in~\cite[Formula 3.39 ff.]{Ge79}. 

\medskip
The proof of the second inequality in \eqref{nutshell} is our main technical result. It rests on the following decomposition for the relative entropy loss densities $g_L(\nu|\mu)$ and $\tilde g_L(\nu|\mu)$. We define the 
\textit{specific entropy loss} by
\begin{equation*}\label{Nup1}
\begin{split}
g^n_L(\nu):=\sum_{\o_{\L_n}\in\{1,\dots,q\}^{\L_n}}\nu(L\one_{\o_{\L_n}})\log\nu(\one_{\o_{\L_n}})
\end{split}
\end{equation*}
and the \textit{specific energy loss} by
\begin{equation*}\label{Nup2}
\begin{split}
\r^n_L(\nu,\mu):=-\sum_{\o_{\L_n}\in\{1,\dots,q\}^{\L_n}}\nu(L\one_{\o_{\L_n}})\log\mu(\one_{\o_{\L_n}}). 
\end{split}
\end{equation*}
Then  
\begin{equation*}\label{Nup3}
\begin{split}
g_L(\nu)=\lim_{n\uparrow\infty}\frac{1}{|\L_n|}g^n_L(\nu)\quad\text{ and }\quad\r_L(\nu,\mu)=\lim_{n\uparrow\infty}\frac{1}{|\L_n|}\r^n_L(\nu,\mu)
\end{split}
\end{equation*}
are their associated densities, whenever their limits exist. Note that, if $\r_L(\nu,\mu)$ and $g_L(\nu)$ are defined, then 
\begin{equation}\label{RelEntLoss}
\begin{split}
g_L(\nu|\mu)=\r_L(\nu,\mu)+g_L(\nu)
\end{split}
\end{equation}
and $g_L(\nu|\mu)$ is given by a limit instead of a limit superior. 
Similar, we define the \textit{approximating specific entropy loss} by
\begin{align*}
\tilde g_L^n(\nu)
=\sum_{i\in\tilde\L_n}\sum_{\D\ni i}\tfrac{1}{|\D\cap\L_n|}\sum_{\xi_{\D}}\int\nu(d\eta)c_\D(\eta_{\L_n},\xi_{\D})\log\tfrac{\nu(\one_{\xi_{\D\cap\L_n}}|{\eta_{\L_n\setminus\D}})}{\nu(\one_{\eta_{\D\cap\L_n}}|{\eta_{\L_n\setminus\D}})},
\end{align*}
the \textit{approximating specific energy loss} by
\begin{align*}
\tilde \r_L^n(\nu|\mu)
=-\sum_{i\in\tilde\L_n}\sum_{\D\ni i}\tfrac{1}{|\D\cap\L_n|}\sum_{\xi_{\D}}\int\nu(d\eta)c_\D(\eta_{\L_n},\xi_{\D})\log\tfrac{\mu(\one_{\xi_{\D\cap\L_n}}|{\eta_{\L_n\setminus\D}})}{\mu(\one_{\eta_{\D\cap\L_n}}|{\eta_{\L_n\setminus\D}})},
\end{align*}
and by
\begin{equation*}\label{Nup3}
\begin{split}
\tilde g_L(\nu)=\lim_{n\uparrow\infty}\frac{1}{|\L_n|}\tilde g^{n}_L(\nu)\quad\text{ and }\quad\tilde \r_L(\nu,\mu)=\lim_{n\uparrow\infty}\frac{1}{|\L_n|}\tilde \r^{n}_L(\nu,\mu)
\end{split}
\end{equation*}
are their associated densities, whenever their limits exist. Again, if $\tilde \r_L(\nu,\mu)$ and $\tilde g_L(\nu)$ are defined, then 
\begin{equation}\label{RelEntLoss}
\begin{split}
\tilde g_L(\nu|\mu)=\tilde \r_L(\nu,\mu)+\tilde g_L(\nu)
\end{split}
\end{equation}
and $\tilde g_L$ is given by a limit instead of a limit superior. 

\medskip
Let us start the analysis of the entropy loss decomposition by giving the following representation result for the specific energy loss.
\begin{prop}\label{SpecificEnergyLoss}
Let $L$ be a well-defined and translation-invariant IPS and $\mu$ a translation-invariant Gibbs measure for the translation-invariant non-null quasilocal specification $\g$. Then, for any $\nu\in\PP_\theta$, 
\begin{equation}\label{ReprEner*}
\begin{split}
\r_L(\nu,\mu)=\tilde \r_L(\nu,\mu)=\sum_{\D\ni o}\sum_{\xi_{\D}}\int\nu(d\eta)c_\D(\eta,\xi_{\D})\frac{1}{|\D|}\log\frac{\g_\D(\one_{\eta_\D}|\eta_{\D^c})}{\g_\D(\one_{\xi_\D}|\eta_{\D^c})}
\end{split}
\end{equation}
and $\nu\mapsto\r_L(\nu,\mu)$ is continuous w.r.t.~the weak topology on $\PP_\theta$.
\end{prop}
We come to our main technical result which states the existence and upper semicontinuity of $\nu\mapsto \tilde g_L(\nu)$. The approach is inspired by the works \cite{Ho71,HiSh75} for the Ising model,  
but much more general. 

\begin{prop}\label{ApproxThmUpperSemi}
Let $L$ be a well-defined IPS generator satisfying conditions~(\ref{tran}), (\ref{fini}), (\ref{cont}), (\ref{mini}) and (\ref{trap}) in Definition~\ref{Gibbs-Attractor-IPS}. Then, $\tilde g_L(\nu)$ exists and is upper semicontinuous on the set of translation-invariant probability measures. 
\end{prop}

In order to have the same statement also for $g_L(\nu)$ we must impose non-nullness.
\begin{prop}\label{ThmUpperSemi}
Let $L$ be as in Proposition~\ref{ApproxThmUpperSemi}. Then, $\tilde g_L(\nu)=g_L(\nu)$ on the set of non-null translation-invariant probability measures.
\end{prop}
We are now in the position to prove our main theorems about irreversible IPS. We start with the version without non-nullness.
\begin{proof}[Proof of Theorem~\ref{MainCorApprox}]
We have the following chain of inequalities,
\begin{equation}\label{UpperSemiApprox}
\begin{split}
0&=\limsup_{k\uparrow\infty}g_L(\nu_{t_k}|\mu)\le\limsup_{k\uparrow\infty}\tilde g_L(\nu_{t_k}|\mu)= \limsup_{k\uparrow\infty}\tilde g_L(\nu_{t_k})+\limsup_{k\uparrow\infty}\tilde \r_L(\nu_{t_k},\mu)\cr
&\le\tilde g_L(\nu_*)+\tilde\r_L(\nu_*,\mu)=\tilde g_L(\nu_*|\mu).
\end{split}
\end{equation}
Here, for the first equality note that the relative entropy density $h(\nu_{t_k}|\mu)$ is non-negative for all $k\in\N$, see for example~\cite{Ge11}. Using Proposition~\ref{EntropyLossLem} and the Condition~\ref{approx zero entropy loss condition} we have that $g_L(\nu_{t_k}|\mu)\le \tilde g_L(\nu_{t_k}|\mu)\le 0$ for all $k\in\N$, and thus the derivative of $h$ given by $g_L$ is non-positive. It hence must approach zero, in fact as a limit, and we do not use non-nullness of $(\nu_{t_k})_{k\in\N}$. The second inequality is justified by Proposition~\ref{EntropyLossLem}, the third and last equalities are by construction. The fourth inequality is a consequence of Proposition~\ref{SpecificEnergyLoss} and Proposition~\ref{ApproxThmUpperSemi}.

By weak closedness of $D_L(\mu)$, under the Condition~\ref{approx zero entropy loss condition}, we then have $\tilde g_L(\nu_*|\mu)\le 0$ and hence from the inequality~\eqref{UpperSemiApprox} it follows that $\tilde g_L(\nu_*|\mu)=0$. Again by Condition~\ref{approx zero entropy loss condition} we then have $\nu_*\in\GG(\g)$, which finishes the proof.
\end{proof}

The version with non-nullness is now a consequence of the fact, that for non-null measures, the approximating relative entropy loss density equals the relative entropy loss density.

\begin{proof}[Proof of Theorem~\ref{MainCor}]
Using Proposition~\ref{ThmUpperSemi}, Proposition~\ref{SpecificEnergyLoss} and Proposition~\ref{EntropyLossLem}  we have that $\tilde g_L(\nu|\mu)=g_L(\nu|\mu)\le 0$ for non-null measures $\nu\in\mathcal P_\theta$. Hence, following the same arguments as in the 
proof of Theorem~\ref{MainCorApprox} we see that $g_L(\nu_*|\mu)=0$ and by Condition~\ref{zero entropy loss condition} we have $\nu_*\in\GG(\g)$.
\end{proof}

The proof of Theorem~\ref{GTilde} is more technical and is given in the following section, together with all other proofs. 

\section{Proofs}\label{Proofs}
In order to increase readability, in this section, we will use the following short-hand notation. Instead of $\nu(\one_{\o_\L})$ we will simply write $\nu(\o_\L)$.
\subsection{Proof of Proposition~\ref{Prop_Non_Null}} 
The essential part of the proof proceeds via a decoupling of a single-site in path space, using
a Girsanov formula to show boundedness of the errors. In this way we
reduce the verification of non-nullness to a single-site situation, where it is easily seen to hold. 
\begin{proof}[Proof of Proposition~\ref{Prop_Non_Null}]
For any $\tau>0$, let us write $\s [0,\tau]$ for the path of a configuration $\s$ in $[0,\tau]$, viewed as a random variable w.r.t.~the law $Q_\nu$ of the Markov process started in $\nu$ and propagated by $(P^L_t)_{t\ge 0}$ in $[0,\tau]$. Further, we denote by $\s(t)$ the projection of $\s[0,\tau]$ to the fixed time $t\in[0,\tau]$ and we write $Q_\o$ instead of $Q_{\d_\o}$.
We will show that there exists an $\d(\tau)>0$, such that for all finite volumes $\L\sm o$, all conditionings $\eta_{\L\sm o}$ and all initial infinite-volume configurations $\o$, there is the lower bound
\begin{equation}\label{GeeFinVol}
\begin{split}
Q_\o(\s_o(\tau)=\eta_o| \s_{\L\sm o}(\tau)=\eta_{\L\sm o})\geq \d(\tau).
\end{split}
\end{equation}
Once we have proved~\eqref{GeeFinVol} then, abbreviating $\nu_t=\nu P_t^L$, we have for all $t\ge \tau$ that
\begin{equation*}
\begin{split}
&\nu_t(\eta_o|\eta_{\L\sm o})=\frac{\int\nu_{t-\tau}(d\o)Q_\o(\s_o(\tau)=\eta_o| \s_{\L\sm o}(\tau)=\eta_{\L\sm o})Q_\o(\s_{\L\sm o}(\tau)=\eta_{\L\sm o})}{\int\nu_{t-\tau}(d\o)Q_\o(\s_{\L\sm o}(\tau)=\eta_{\L\sm o})}\geq \d(\tau)
\end{split}
\end{equation*}
and the almost-sure lower bound carries over by martingale convergence under the limit $\L\uparrow \Z^d$, i.e., 
\begin{equation}\label{GeeFinVolNN}
\begin{split}
&\nu_t(\eta_o| \eta_{o^c})
=   \lim_{\L\uparrow\Z^d}\nu_t(\eta_o| \eta_{\L\sm o})\geq \d(\tau).
\end{split}
\end{equation}
To prove~\eqref{GeeFinVol}, we introduce a reference process with generator 
$L^b$ of the form
$$L^b f(\eta)=\sum_{i\in\Z^d}\sum_{\xi_i}b_i(\eta,\xi_i)[f(\xi_i\eta_{i^c})-f(\eta)].$$
$L^b$ will be a non-translation-invariant finite-volume perturbation of $L$ obeying all the other conditions we assumed for $L$. 
By this we mean more precisely, that the rates are equal, $b_{i}=c_{i}$, {\em except} (possibly) for finitely many flipping sites $i\in\D$. 
We also assume that the reachability functions $d$ are the same for both generators 
 $L$ and $L^b$.
 
Denote by $Q_\o^b$ the measure on path space of the process with generator $L^b$ started in $\o$. 
Then (compare~\cite{MaReVe01} in the context of a discussion of fluctuation laws or~\cite[Theorem 4]{Da93} in the context of large deviations for IPS) 
we have for the Radon-Nikodym derivative in the space of c{\'a}dl{\'a}g 
paths the useful formula  
\begin{equation}\label{Gee0}
\begin{split}
&\frac{dQ_{\o}}{dQ^b_{\o}}(\s[0,\tau])=\exp\Bigl( 
-\int_{0}^\tau \lambda(\s(s))ds+\sum_{{s\in [0,\tau]}\atop{\s_\D(s^-)\neq \s_\D(s)}}\sum_{i\in\D}\log 
\frac{c_i(\s(s^-),\s_i(s))}{b_i(\s(s^-),\s_i(s))}
\Bigr).
\end{split}
\end{equation}
In the first term we have written 
$$\lambda(\eta):=\sum_{i\in\D}\big(c_i(\eta)-b_i(\eta)\big)$$ 
for the difference in the total rates of exiting a configuration $\eta$ taken for the 
different generators. This sum involves only finitely many updating sets $i\in\D$, as the rates of the two generators differ only in a finite volume. Moreover, the integral is well-defined almost surely since almost surely, in a finite volume, there are only finitely many updates. In the second term on the r.h.s.~of~\eqref{Gee0}, there is a sum over jumping times of the process which again has to be taken only in finite volume $\D$, as the rates coincide except in a finite volume and again, in finite volume, almost surely, there are only finitely many updates. 
Hence, the Radon-Nikodym derivative is a local function in path-space. 
Spelling out the finite-volume conditional probability of the time-evolved measure on the l.h.s.~
of~\eqref{GeeFinVol}, we would like the bring the reference process for a suitable generator $L^b$ into play. Let us therefore abbreviate both parts of the Radon-Nikodym derivative as 
\begin{equation}\label{Gee2}
\begin{split}
&a(\s[0,\tau]):=\exp\Bigl( -\int_{0}^\tau \lambda(\s(s))ds\Bigr)\cr
&A(\s[0,\tau]):=\exp\Bigl( \sum_{{s\in [0,\tau]}\atop{\s_\D(s^-)\neq \s_\D(s)}}\sum_{i\in\D}\log 
\frac{c_i(\s(s^-),\s_i(s))}{b_i(\s(s^-),\s_i(s))}
\Bigr).
\end{split}
\end{equation}
Then, we have 
\begin{equation*}
\begin{split}
Q_\o(\s_o(\tau)=\eta_o| \s_{\L\sm o}(\tau)=\eta_{\L\sm o})&=\frac{
Q_\o(\one_{\s_o(\tau)=\eta_o}\one_{ \s_{\L\sm o}(\tau)=\eta_{\L\sm o}})
}{Q_\o(\one_{ \s_{\L\sm o}(\tau)=\eta_{\L\sm o}})}\cr
&=\frac{
Q_\o^b(a(\s[0,\tau])A(\s[0,\tau])\one_{\s_o(\tau)=\eta_o}\one_{ \s_{\L\sm o}(\tau)=\eta_{\L\sm o}})
}{Q_\o^b(a(\s[0,\tau])A(\s[0,\tau])\one_{ \s_{\L\sm o}(\tau)=\eta_{\L\sm o}})}.
\end{split}
\end{equation*}
From now on we assume that $Q$ denotes the path-measure with an arbitrary 
fixed initial configuration $\o$ which we will drop now in the notation. 
By the boundedness of the rates, both of $L$ and the reference process $L^b$, we have 
that for each $\tau$, there exists an $\k(\tau)>0$ such that for almost all paths $\s[0,\tau]$ there is 
a deterministic upper and lower bound
$\k(\tau)\leq a(\s[0,\tau])\leq \frac{1}{\k(\tau)}$. 
So we can remove this factor and obtain the lower bound 
\begin{equation*}
\begin{split}
&Q(\s_o(\tau)=\eta_o| \s_{\L\sm o}(\tau)=\eta_{\L\sm o})\geq \k^2(\tau)\frac{
Q^b(A(\s[0,\tau])\one_{\s_o(\tau)=\eta_o}\one_{ \s_{\L\sm o}(\tau)=\eta_{\L\sm o}})
}{Q^b(A(\s[0,\tau])\one_{ \s_{\L\sm o}(\tau)=\eta_{\L\sm o}})}.
\end{split}
\end{equation*}
We would like to reduce the problem to a single-site problem, and we have to pay 
 in terms of the appearance of the more dangerous part  $A(\cdot)$. It might look 
 problematic, as it is not bounded uniformly in all 
possible trajectories, as the number of jumps in the finite volume $\D$ can be arbitrarily large with positive probability.   
Nevertheless, we can control it in terms of tails of Poisson variables with bounded intensities, as we will see. 

To carry this out explicitly we have to choose first the decoupling generator $L^b$ in such a way, that a modification is made for all update sites $i$ which might be influencing the state of the process at site $o$. But also, conversely, might be influenced by the state of the process at site $o$. To do so, denote the maximal rate of $L$ by 
$\hat c:=\sup_{\eta,\xi_{i}}c_i(\eta,\xi_{i})$ 
and put 
\begin{equation*}
\begin{split}
&b_{i}(\eta,\xi_{i}):=\hat c d (\eta_{i},\xi_{i})\one_{\D}(i)+ c_{i}(\eta,\xi_{i})\one_{\D^c}(i)
\end{split}
\end{equation*}
The main point of this definition is that the processes $\s_o[0,\tau]$ 
and $\s_{o^c}[0,\tau]$ are now independent under the joint law $Q^b$, for any fixed starting configuration. 
Let us come to the treatment of the unbounded term $A$. 
By our choice we have from $b_{i}\geq c_{i}$ the deterministic upper bound $A\leq 1$. 
Making use of the minimal transition rate property
$\inf_{\eta,\xi_{i}:\, c_{i}(\eta,\xi_{i})>0}c_{i}(\eta,\xi_{i})>0$, one of our general assumptions which is automatic by finite range, we also have the lower bound in terms of the total number of jumps around $0$, i.e., 
\begin{equation*}
\begin{split}
&A(\s[0,\tau])\geq e^{ -R N_{\D}(\tau)}. 
\end{split}
\end{equation*}
Here, for a general volume $W$ we write  
$N_{W}(\tau):=\#\{s\in [0,\tau], \o_{W}(s^-)\neq \o_{W}(s)\}$
for the number of updates in $W$ until time $T$. 
Here, $R>0$ is an obvious finite constant related to the minimal and maximal rate.    
To summarize, what we have obtained so far is 
\begin{equation*}
\begin{split}
&Q(\s_o(\tau)=\eta_o| \s_{\L\sm o}(\tau)=\eta_{\L\sm o})
\geq \k^2(\tau)\frac{
Q^b(e^{ -R N_{\D}(\tau)}\one_{\s_o(\tau)=\eta_o}\one_{\s_{\L\sm o}(\tau)=\eta_{\L\sm o}})
}{Q^b(\one_{ \s_{\L\sm o}(\tau)=\eta_{\L\sm o}})}.
\end{split}
\end{equation*}
Using that $N_{\D}=N_{o}+N_{\D\ba o}$ and the independence 
of the processes at site $o$ and away from site $o$ under  $Q^b$, 
the term in question factorizes and we have 
\begin{equation}\label{Gee19}
\begin{split}
&\frac{
Q^b(e^{ -R N_{\D}(\tau)}\one_{\s_o(\tau)=\eta_o}\one_{ \s_{\L\sm o}(\tau)=\eta_{\L\sm o}})
}{Q^b(\one_{ \s_{\L\sm o}(\tau)=\eta_{\L\sm o}})}\cr
&=Q^b(e^{ -R N_o(\tau)}\one_{\s_o(\tau)=\eta_o})\times 
Q^b(e^{ -R N_{\D\ba o}(\tau)} | \s_{\L\sm o}(\tau)=\eta_{\L\sm o}).
\end{split}
\end{equation}
The first term is just an expression in the irreducible time-homogeneous single-site Markov chain 
with transition rates $\hat c d$.  
To be explicit, for any single-site initial condition $\o_o$ we write lower bounds
\begin{equation}\label{Gee17}
\begin{split}
Q^b(e^{ -R N_o(\tau)}\one_{\s_o(\tau)=\eta_o})
&\geq e^{- R n}Q^b(N_o(\tau)\leq n, \s_o(\tau)=\eta_o)\cr
&\geq e^{- R n}\bigl(Q^b(\s_o(\tau)=\eta_o)-Q^b(N_o(\tau)> n)\bigr).
\end{split}
\end{equation}
First, by irreducibility, for any $\tau>0$ and all $\eta_o$, there is a lower 
bound  $Q^b(\s_o(\tau)=\eta_o)\geq \tilde\r(\tau)>0$. Next, 
$Q^b(N_o(\tau)> n)$ is bounded from above by the tail of a Poisson distribution, 
and hence, now choosing $n$ sufficiently large but finite, 
gives a strictly positive uniform lower bound $\r(\tau)$ on the r.h.s.~of~\eqref{Gee17}.

\medskip
The second term of~\eqref{Gee19} is bounded below by
\begin{equation*}
\begin{split}
e^{- R n}(1-Q^b( N_{\D\ba o}(\tau)> n &| \s_{\L\sm o}(\tau)=\eta_{\L\sm o}))\cr
&\geq e^{- R n}(1-\sum_{j\in \D\ba o}Q^b( N_{j}(\tau)> \frac{n}{|\D\ba o|} | \s_{\L\sm o}(\tau)=\eta_{\L\sm o})).
\end{split}
\end{equation*}
Considering one summand of the finitely many sites 
$j\in \D\ba o$, we condition on the behavior of the path away from $j$ and write  
\begin{equation}\label{GeeFinal}
\begin{split}
Q^b(& N_{j}(\tau)> m | \s_{\L\sm o}(\tau)=\eta_{\L\sm o})\cr
&=\int Q^b( d\o_{j^c}[0,\tau] | \s_{\L\sm o}(\tau)=\eta_{\L\sm o})\times
Q^b( N_{j}(\tau)> m | \s_{j^c}[0,\tau]=\o_{j^c}[0,\tau]).
\cr
\end{split}
\end{equation}
Now, for any fixed realization of the path $\o_{j^c}[0,\tau]$ away from $j$ 
the term under the integral on the r.h.s.~is the tail of the counting variable $N_{j}(\tau)$ 
of a single-site Markov chain at site $j$ with time-inhomogeneous rates which 
are given by the behavior of the path $\o_{j^c}[0,\tau]$ in a finite neighborhood around the site $j$. 
These rates are not explicit, but all that matters is that they are uniformly bounded by  
$\hat c\tau$. 
Hence, the r.h.s.~of~\eqref{GeeFinal} is dominated from above by the corresponding tail 
of a Poisson variable 
with parameter $\hat c\tau$.
Now choose $m=n/|\D\ba o|$ sufficiently large but finite, to finish the 
proof of non-nullness. 

\medskip
Finally, the non-nullness carries over to any subsequential limiting measure since the bound $\d(\tau)$ is uniform in time, volume and the measure.
\end{proof}

\subsection{Proof of Proposition~\ref{EntropyLossLem}}
The proof is based on a finite-volume argument for an approximating dynamics which depends on the infinite-volume time-stationary measure $\mu$, using Jensen's inequality. We note that this proof would become
much simpler in the reversible setting, where the reversible measure and the rates can be related configuration wise via the detailed balance equations, see the proof of Theorem~\ref{GTilde}. In the irreversible case, we loose this identification and the arguments become more involved. 

Before we enter the proof, let us give a representation of the relative entropy loss $g_L^n(\nu|\mu)$ that we will use for calculations in what follows. 
\begin{lem}\label{Lem_Rep}
Let $c^\nu_\D(\eta_{\L_n},\xi_{\D\cap\L_n}):=\sum_{\xi_{\D\sm\L_n}}\int\nu(d\s|\eta_{\L_n})c_\D(\eta_{\L_n}\s_{\L_n^c},\xi_{\D\cap\L_n}\xi_{\D\sm\L_n})$, then we have that  
\begin{align*}
g_L^n(\nu|\mu)=\sum_{i\in\L_n}\sum_{\D\ni i}\tfrac{1}{|\D\cap\L_n|}\sum_{\xi_{\D\cap\L_n}}\sum_{\eta_{\L_n}}\nu({\eta_{\L_n}})c^\nu_\D(\eta_{\L_n},\xi_{\D\cap\L_n})\log\tfrac{\nu({\xi_{\D\cap\L_n}}|{\eta_{\L_n\setminus\D}})\mu({\eta_{\D\cap\L_n}}|{\eta_{\L_n\setminus\D}})}{\nu({\eta_{\D\cap\L_n}}|{\eta_{\L_n\setminus\D}})\mu({\xi_{\D\cap\L_n}}|{\eta_{\L_n\setminus\D}})}.
\end{align*}
\end{lem}
\begin{proof}[Proof of Lemma~\ref{Lem_Rep}]
We prove by direct calculation. We have that
\begin{equation*}\label{EntropyCompute_1}
\begin{split}
g_L^n(\nu,\mu)&=\sum_{\o_{\L_n}}\nu(L\one_{\o_{\L_n}})\log\frac{\nu(\o_{\L_n})}{\mu(\o_{\L_n})}\cr
&=\sum_{\o_{\L_n}}\log\frac{\nu(\o_{\L_n})}{\mu(\o_{\L_n})}\sum_{\D:\, \D\cap\L_n\neq\emptyset}\sum_{\xi_\D}\int\nu(d\eta)c_\D(\eta,\xi_\D)[\one_{\o_{\L_n}}(\xi_\D\eta_{\D^c})-\one_{\o_{\L_n}}(\eta)]\cr
&=\sum_{\D:\, \D\cap\L_n\neq\emptyset}\sum_{\xi_\D}\int\nu(d\eta)c_\D(\eta,\xi_\D)\log\frac{\nu(\xi_{\D\cap\L_n}\eta_{\L_n\sm\D})\mu(\eta_{\L_n})}{\nu(\eta_{\L_n})\mu(\xi_{\D\cap\L_n}\eta_{\L_n\sm\D})}\cr
&=\sum_{i\in\L_n}\sum_{\D\ni i}\tfrac{1}{|\D\cap\L_n|}\sum_{\xi_\D}\int\nu(d\eta)c_\D(\eta,\xi_\D)\log\frac{\nu({\xi_{\D\cap\L_n}}|{\eta_{\L_n\setminus\D}})\mu({\eta_{\D\cap\L_n}}|{\eta_{\L_n\setminus\D}})}{\nu({\eta_{\D\cap\L_n}}|{\eta_{\L_n\setminus\D}})\mu({\xi_{\D\cap\L_n}}|{\eta_{\L_n\setminus\D}})}.
\end{split}
\end{equation*}
Replacing the rates, implies the desired representation.
\end{proof}
Also we will often need estimates on the logarithmic terms in $g_L^n(\nu|\mu)$. For this the nun-nullness is a sufficient assumption. We have the following estimate. 
\begin{lem}\label{Lem_Condition}
Let $\D\subset\L\Subset\Z^d$ and $\nu\in\mathcal P_\theta$ be non-null with parameter $\d>0$, then for all $\eta,\xi\in\O$ we have that
\begin{align*}
|\log\tfrac{\nu({\xi_{\D}}|{\eta_{\L\setminus\D}})}{\nu({\eta_{\D}}|{\eta_{\L\setminus\D}})}|\le |\D|\log\frac 1 \d.
\end{align*}
In particular, for $\mu\in\GG(\g)$ with non-null specification $\g$, the same estimate holds. 
\end{lem}
\begin{proof}[Proof of Lemma~\ref{Lem_Condition}]
Let $i_1,\dots,i_k$ be any numbering of the sites in $\D$, and denote $[i_j,i_k]=\{i_j,i_{j+1},\dots, i_k\}$, then 
by the chain rule of conditional probabilities
\begin{align*}
\nu({\eta_{\D}}|{\eta_{\L\setminus\D}})=\prod_{j=1,\dots, k-1}\nu({\eta_{i_j}}|\eta_{[i_{j+1},i_k]}{\eta_{\L\setminus\D}}).
\end{align*}
Now, for every $i_j$ we can expand the conditioning in an elementary way and write
\begin{equation}\label{NonNullEst}
\begin{split}
&\nu({\eta_{i_j}}|\eta_{[i_{j+1},i_k]}{\eta_{\L\setminus\D}})=\frac{\int\nu(d\s)\nu({\eta_{i_j}}\eta_{[i_{j+1},i_k]}{\eta_{\L\setminus\D}}|\s_{\L^c\cup[i_1,i_{j-1}]})}{\int\nu(d\s)\nu(\eta_{[i_{j+1},i_k]}{\eta_{\L\setminus\D}}|\s_{\L^c\cup[i_1,i_{j-1}]})}\cr
&=\frac{\int\nu(d\s)\frac{\nu({\eta_{i_j}}\eta_{[i_{j+1},i_k]}{\eta_{\L\setminus\D}}|\s_{\L^c\cup[i_1,i_{j-1}]})}{\nu(\eta_{[i_{j+1},i_k]}{\eta_{\L\setminus\D}}|\s_{\L^c\cup[i_1,i_{j-1}]})}\nu(\eta_{[i_{j+1},i_k]}{\eta_{\L\setminus\D}}|\s_{\L^c\cup[i_1,i_{j-1}]})}{\int\nu(d\s)\nu(\eta_{[i_{j+1},i_k]}{\eta_{\L\setminus\D}}|\s_{\L^c\cup[i_1,i_{j-1}]})}\cr
&=\frac{\int\nu(d\s)\nu({\eta_{i_j}}|\eta_{[i_{j+1},i_k]}{\eta_{\L\setminus\D}}\s_{\L^c\cup[i_1,i_{j-1}]})\nu(\eta_{[i_{j+1},i_k]}{\eta_{\L\setminus\D}}|\s_{\L^c\cup[i_1,i_{j-1}]})}{\int\nu(d\s)\nu(\eta_{[i_{j+1},i_k]}{\eta_{\L\setminus\D}}|\s_{\L^c\cup[i_1,i_{j-1}]})}\ge \d
\end{split}
\end{equation}
which gives the desired bound. If $\nu$ is replaced by $\mu\in\GG(\g)$ for a non-null specification $\g$, then, in equation~\eqref{NonNullEst}, using the DLR equation, the term $\nu({\eta_{i_j}}|\eta_{[i_{j+1},i_k]}{\eta_{\L\setminus\D}}\s_{\L^c\cup[i_1,i_{j-1}]})$ can be replaced by $\g({\eta_{i_j}}|\eta_{[i_{j+1},i_k]}{\eta_{\L\setminus\D}}\s_{\L^c\cup[i_1,i_{j-1}]})$ which leads to the same bound. 
\end{proof}
\begin{remark}
Let us mention that throughout the manuscript, the non-nullness condition is stronger than necessary. What is really needed is that
\begin{equation*}
\begin{split}
\sup_{\L_n\ni o}\sum_{\D\ni o}\sum_{\xi_\D}\int\nu(d\eta)c^{\xi_\D}_\D(\eta)\log\frac{\nu(\eta_{\D\cap\L_n}|\eta_{\L_n\setminus\D})}{\nu(\xi_{\D\cap\L_n}|\eta_{\L_n\setminus\D})}<\infty
\end{split}
\end{equation*}
which is implied if $\nu$ is non-null.  
\end{remark}

\begin{proof}[Proof of Proposition~\ref{EntropyLossLem}]
Let us first show that $g_L(\nu|\mu)\le\tilde g_L(\nu|\mu)$ for all $\nu\in\mathcal P_\theta$.  Indeed, using $\log x\le x$ and Lemma~\ref{Lem_Condition}, we can estimate 
\begin{equation}\label{TildeLarger}
\begin{split}
g_L^n(\nu|\mu)-\tilde g_L^n(\nu|\mu)&=\sum_{i\in\L_n\sm\tilde\L_n}\sum_{\D\ni i}\tfrac{1}{|\D\cap\L_n|}\sum_{\xi_{\D\cap\L_n}}\sum_{\eta_{\L_n}}\nu( {\eta_{\L_n}})c^\nu_\D(\eta_{\L_n},\xi_{\D\cap\L_n})\cr
&\hspace{4cm}\times\log\tfrac{\nu( {\xi_{\D\cap\L_n}}| {\eta_{\L_n\setminus\D}})\mu( {\eta_{\D\cap\L_n}}| {\eta_{\L_n\setminus\D}})}{\nu( {\eta_{\D\cap\L_n}}| {\eta_{\L_n\setminus\D}})\mu( {\xi_{\D\cap\L_n}}| {\eta_{\L_n\setminus\D}})}\cr
&\le\sum_{i\in\L_n\sm\tilde\L_n}\sum_{\D\ni i}\tfrac{1}{|\D\cap\L_n|}\sum_{\xi_{\D\cap\L_n}}\sum_{\eta_{\L_n}}\nu( {\eta_{\L_n}})c^\nu_\D(\eta_{\L_n},\xi_{\D\cap\L_n})\cr
&\hspace{4cm}\times\big[\tfrac{\nu( {\xi_{\D\cap\L_n}}| {\eta_{\L_n\setminus\D}})}{\nu( {\eta_{\D\cap\L_n}}| {\eta_{\L_n\setminus\D}})}+|\D\cap\L_n|\log\frac{1}{\d} \big]\cr
&=\sum_{i\in\L_n\sm\tilde\L_n}\sum_{\D\ni i}\tfrac{1}{|\D\cap\L_n|}\sum_{\xi_{\D\cap\L_n}}\sum_{\eta_{\L_n}}c^\nu_\D(\eta_{\L_n},\xi_{\D\cap\L_n})\cr
&\hspace{3cm}\times\big[\nu( {\xi_{\D\cap\L_n}}{\eta_{\L_n\setminus\D}})+\nu( {\eta_{\L_n}})|\D\cap\L_n|\log\frac{1}{\d} \big]\cr
&= |\L_n\sm\tilde\L_n|\sum_{\D\ni o}c_\D(q^{|\D|}+ \log\frac{1}{\d})
\end{split}
\end{equation}
where we also used translation invariance in the last equality. Since we assume only finitely many types of transitions, by the definition of $\L_n$ and $\tilde\L_n$, this tends to zero in the density limit.

\medskip
Now we show $g_L(\nu|\mu)\le 0$. For this, consider the approximating finite-volume generator
\begin{equation*}\label{IPSappr}
L_n\one_{\o_{\L_n}}(\eta_{\L_n})=\sum_{\D:\,\D\subset{\L_n}}\sum_{\xi_\D}c^n_\D(\eta_{\L_n},\xi_\D)[\one_{\o_{\L_n}}(\xi_{\D}\eta_{{\L_n}\setminus\D})-\one_{\o_{\L_n}}(\eta_{\L_n})]
\end{equation*}
where for $\D\subset{\L_n}$ the approximating rates are defined by
$$c_{\D}^n(\eta_{\L_n},\xi_{\D})=\sum_{\D'\Subset\L_n^c}\int\mu(d\s_{\L_n^c}|\eta_{\L_n})\sum_{\zeta_{\D'}}c_{\D\cup\D'}(\eta_{\L_n}\s_{\L_n^c},\xi_\D\zeta_{\D'})$$ 
with $\mu$ the time-stationary Gibbs measure for the infinite-volume generator $L$. Note that $L_n$ is a well-defined finite-volume generator. The finite-volume rates are obtained from the infinite-volume rates in terms of the following two averaging operations. First, a conditional average over the part of the configuration outside of the finite volume $\L_n$ given the part $\eta_{\L_n}$ inside $\L_n$, w.r.t.~the invariant measure $\mu$ in the infinite volume. Second, a sum over 
the part $\D'$ of the update-set $\D\cup \D'$ and the corresponding configurations $\zeta_{\D'}$ outside of $\Lambda$.  
For example if $L$ is the generator of the exclusion process, then the particle number in $\L$ is not preserved by $L_n$ and particles are produced and vanish at the boundary.

This construction in particular implies that $\mu$, as a measure on $\{1,\dots,q\}^{\L_n}$, is invariant w.r.t.~$L_n$. Indeed, we can write
\begin{equation*}\label{IPSappr2}
\begin{split}
0=\mu(L\one_{\o_\L})=\sum_{\eta_{\L_n}}\mu(\eta_{\L_n})\int\mu (d\s_{\L_n^c}|\eta_{\L_n} ) (L \one_{\o_{\L_n}})(\eta_{\L_n}
\s_{\L^c})
\end{split}
\end{equation*}
where the first equality is the time-stationarity of $\mu$ w.r.t.~the infinite-volume dynamics and
the second equality is the decomposition of the infinite-volume measure into its finite-volume marginal
and its conditional probability outside. Hence it suffices to show that $L_n\one_{\o_{\L_n}}(\eta_{\L_n})=\int\mu (d\s_{\L_n^c}|\eta_{\L_n} ) (L\one_{\o_{\L_n}})(\eta_{\L_n}\s_{\L_n^c})$. But this is true since 
\begin{equation*}\label{IPSappr2}
\begin{split}
&L_n\one_{\o_{\L_n}}(\eta_{\L_n})=\sum_{\D:\,\D\subset\L_n}\sum_{\xi_\D}c^n_\D(\eta_{\L_n},\xi_\D)[\one_{\o_{\L_n}}(\xi_\D\eta_{\L_n\setminus\D})-\one_{\o_{\L_n}}(\eta_{\L_n})]\cr
&=\int\mu(d\s_{\L_n^c}|\eta_{\L_n})\sum_{\D:\,\D\subset\L_n}\sum_{\D'\Subset\L_n^c}\sum_{\xi_\D}\sum_{\zeta_{\D'}}c_{\D\cup\D'}(\eta_{\L_n}\s_{\L_n^c},\xi_\D\zeta_{\D'})[\one_{\o_{\L_n}}(\xi_\D\eta_{\L_n\setminus\D})-\one_{\o_{\L_n}}(\eta_{\L_n})]\cr
&=\int\mu(d\s_{\L_n^c}|\eta_{\L_n})\sum_{\D'':\, \D''\cap\L_n\neq\emptyset}\sum_{\xi_{\D''}}c_\D(\eta_{\L_n}\s_{\L_n^c},\xi_{\D''})[\one_{\o_{\L_n}}(\xi_{\D''\cap\L_n}\eta_{\L_n\setminus\D''})-\one_{\o_{\L_n}}(\eta_{\L_n})]\cr
&=\int\mu(d\s_{\L_n^c}|\eta_{\L_n})L\one_{\o_{\L_n}}(\eta_{\L_n}\s_{\L_n^c})
\end{split}
\end{equation*}
where we have used, that every $\D''$ appearing in the sum in the third line has a unique decomposition into $\D\subset\L_n$ and $\D'\subset\L_n^c$ appearing in the second line.

\medskip
Now, let $(P_t^{n})_{t\ge0}$ denote the semigroup associated to $L_n$, then by Jensen's inequality applied to the concave function $\Psi(u):=-u\log u+u-1$ we have 
\begin{equation*}\label{Relative__Entropy_Approx_Jensen}
\begin{split}
h_{\L_n}(\nu P_t^{n}|\mu)&=-\sum_{\eta_{\L_n}}\mu(\eta_{\L_n})\Psi(\frac{\nu P_t^{n}(\eta_{\L_n})}{\mu(\eta_{\L_n})})
\le-\sum_{\eta_{\L_n}}\mu(\eta_{\L_n})\Psi(\frac{\nu(\eta_{\L_n})}{\mu(\eta_{\L_n})})=h_{\L_n}(\nu|\mu).
\end{split}
\end{equation*} 
This is a standard argument for finite Markov processes, see for example \cite[Theorem 3.A3]{Ge11}. Consequently,  the derivative 
\begin{equation*}\label{Relative__Entropy_Approx_Jensen}
\begin{split}
\frac{d}{dt}_{|t=0}h_{\L_n}(\nu P_t^{n}|\mu)=\sum_{\D:\,\D\subset\L_n}\sum_{\xi_\D}\sum_{\eta_{\L_n}}\nu(\eta_{\L_n})c^n_\D(\eta_{\L_n},\xi_\D)\log\frac{\nu(\eta_{\L_n})\mu(\xi_{\D}\eta_{{\L_n}\setminus\D})}{\nu(\xi_{\D}\eta_{\L_n\setminus\D})\mu(\eta_{\L_n})}
\end{split}
\end{equation*}
must be non-positive. What remains to show, is that the error produced by the approximation of the dynamics is of boundary order. 
Recall the short-hand notation for the rates from Lemma~\ref{Lem_Rep}. Then we have the following estimate,
\begin{equation*}\label{Relative__Entropy_Approx_Jensen_1}
\begin{split}
|\frac{d}{dt}_{|t=0}h_{\L_n}(\nu P_t^{L}|\mu)-&\frac{d}{dt}_{|t=0}h_{\L_n}(\nu P_t^{n}|\mu)|\le\sum_{i\in\L_n}\sum_{\D\ni i}\tfrac{1}{|\D\cap\L_n|}\sum_{\xi_{\D\cap\L_n}}\sum_{\eta_{\L_n}}\nu(\eta_{\L_n})\cr
&\times|c^\nu_\D(\eta_{\L_n},\xi_{\D\cap\L_n})-c^\mu_\D(\eta_{\L_n},\xi_{\D\cap\L_n})||\log\tfrac{\nu(\xi_{\D\cap\L_n}\eta_{\L_n\setminus\D})\mu(\eta_{\L_n})}{\nu(\eta_{\L_n})\mu(\xi_{\D\cap\L_n}\eta_{\L_n\setminus\D})}|
\end{split}
\end{equation*}
where we resolved the rates $c^n$ into the notation $c^\mu$. Roughly speaking, by the continuity of the rates, the distance between the rates becomes small for updates in the bulk, and boundary terms can be estimated using non-nullness. To make this precise, we split the first sum into a sum $\sum_{i\in\tilde\L_n}$ of bulk terms and a sum $\sum_{i\in\L_n\sm\tilde\L_n}$ of boundary terms. The sum of boundary terms, using Lemma~\ref{Lem_Condition}, can be bounded from above by 
\begin{equation}\label{Relative__Entropy_Approx_Jensen_2}
\begin{split}
|\L_n\sm\tilde\L_n|4\log \frac 1 \d \sum_{\D\ni o}c_\D
\end{split}
\end{equation}
which tends to zero in the density limit. As for the sum of bulk terms, let $n$ be sufficiently large, such that for all $\D$ with $\D\cap\tilde\L_n\neq\emptyset$ and $c_\D>0$ we have $\D\subset\L_n$. Then, again using Lemma~\ref{Lem_Condition}, we can estimate from above by 
\begin{equation*}\label{Relative__Entropy_Approx_Jensen_3}
\begin{split}
&2\log\frac 1 \d |\tilde\L_n|\sum_{\D\ni o:\, c_\D>0}q^{|\D|}\sup_{\xi,\eta,\s,\s'}|c_\D(\eta_{\L_n}\s_{\L_n^c},\xi_{\D})-c_\D(\eta_{\L_n}\s'_{\L_n^c},\xi_{\D})|. 
\end{split}
\end{equation*}
But since we assumed the rates to be uniformly continuous, this term tends to zero in the density limit as $n$ tends to infinity. This finishes the proof. 
\end{proof}

\subsection{Proof of Proposition~\ref{SpecificEnergyLoss}}
The proof is based on a proper separation of bulk and boundary terms. The main argument then rests on the continuity of the rates and translation invariance. Note, that in this proof, we do not require only finitely many types of transitions. 
\begin{proof}[Proof of Proposition~\ref{SpecificEnergyLoss}]
Considering the proof of Lemma~\ref{Lem_Rep}, we have the following representation of the finite-volume specific energy loss 
\begin{equation*}\label{EntropyCompute}
\begin{split}
\r_L^n&(\nu,\mu)=\sum_{i\in\L_n}\sum_{\D\ni i}\tfrac{1}{|\D\cap\L_n|}\sum_{\xi_{\D\cap\L_n}}\sum_{\eta_{\L_n}}\nu({\eta_{\L_n}})c^\nu_\D(\eta_{\L_n},\xi_{\D\cap\L_n})\log\tfrac{\mu({\eta_{\D\cap\L_n}}|{\eta_{\L_n\setminus\D}})}{\mu({\xi_{\D\cap\L_n}}|{\eta_{\L_n\setminus\D}})}.
\end{split}
\end{equation*}
Using the exact same arguments as in the estimate for the boundary term in~\eqref{Relative__Entropy_Approx_Jensen_2}, we see that $|\r_L^n(\nu,\mu)-\tilde\r_L^n(\nu,\mu)|=o(|\L_n|)$ which proves that $\r_L(\nu,\mu)=\tilde \r_L(\nu,\mu)$ if the limit exists. In order to see that the limit exists and has the desired form, note that, by translation-invariance, the representation on the r.h.s.~of \eqref{ReprEner*} can be written as
\begin{equation*}\label{EntropyChange2}
\begin{split}
\frac{1}{|\L_n|}\sum_{i\in\L_n}\sum_{\D\ni i}\frac{1}{|\D|} \sum_{\xi_\D}\int\nu(d\eta)c_\D(\eta,\xi_{\D})\log\frac{\g_\D(\eta_\D|\eta_{\D^c})}{\g_\D(\xi_\D|\eta_{\D^c})}=:R_L(\nu,\mu).\cr
\end{split}
\end{equation*}
Thus, resolving the short-hand notation for the rates $c_\D^\nu$, introduced in Lemma~\ref{Lem_Rep}, the non-normalized finite-volume difference can be expressed as
\begin{equation}\label{DifferenceEnergy}
\begin{split}
&|\L_n|R_L(\nu,\mu)-\r_L^n(\nu,\mu)\cr
&=\sum_{i\in\L_n}\sum_{\D\ni i}\sum_{\xi_\D}\frac{1}{|\D|}\int\nu(d\eta)c_\D(\eta,\xi_{\D})\log\frac{\g_\D(\eta_\D|\eta_{\D^c})\mu(\xi_{\D\cap\L_n}|\eta_{\L_n\setminus\D})}{\g_\D(\xi_\D|\eta_{\D^c})\mu(\eta_{\D\cap\L_n}|\eta_{\L_n\setminus\D})}\cr
&\quad+\sum_{i\in\L_n}\sum_{\D\ni i}\sum_{\xi_\D}\int\nu(d\eta)(\frac{1}{|\D|}-\frac{1}{|\D\cap\L_n|})c_\D(\eta,\xi_\D)\log\frac{\mu(\eta_{\D\cap\L_n}|\eta_{\L_n\setminus\D})}{\mu(\xi_{\D\cap\L_n}|\eta_{\L_n\setminus\D})}
\end{split}
\end{equation}
and it suffices to show that this difference is of the order $o(|\L_n|)$. 
We would like to separate boundary terms from bulk terms. For this, let us fix a finite set $o\in \G\Subset\Z^d$ and rewrite~\eqref{DifferenceEnergy} as
\begin{equation*}
\begin{split}
&\sum_{i\in\L_n:\, \G+i\subset\L_n}\sum_{\D\ni i:\, \D\not\subset\G+i}\sum_{\xi_\D}\frac{1}{|\D|}\int\nu(d\eta)c_\D(\eta,\xi_{\D})\log\frac{\g_\D(\eta_\D|\eta_{\D^c})\mu(\xi_{\D\cap\L_n}|\eta_{\L_n\setminus\D})}{\g_\D(\xi_\D|\eta_{\D^c})\mu(\eta_{\D\cap\L_n}|\eta_{\L_n\setminus\D})}\cr
&\quad+\sum_{i\in\L_n:\, \G+i\subset\L_n}\sum_{\D\ni i:\, \D\subset\G+i}\sum_{\xi_\D}\frac{1}{|\D|}\int\nu(d\eta)c_\D(\eta,\xi_{\D})\log\frac{\g_\D(\eta_\D|\eta_{\D^c})\mu(\xi_{\D}|\eta_{\L_n\setminus\D})}{\g_\D(\xi_\D|\eta_{\D^c})\mu(\eta_{\D}|\eta_{\L_n\setminus\D})}\cr
&\quad+\sum_{i\in\L_n:\, \G+i\not\subset\L_n}\sum_{\D\ni i}\sum_{\xi_\D}\frac{1}{|\D|}\int\nu(d\eta)c_\D(\eta,\xi_{\D})\log\frac{\g_\D(\eta_\D|\eta_{\D^c})\mu(\xi_{\D\cap\L_n}|\eta_{\L_n\setminus\D})}{\g_\D(\xi_\D|\eta_{\D^c})\mu(\eta_{\D\cap\L_n}|\eta_{\L_n\setminus\D})}\cr
&\quad+\sum_{i\in\L_n:\, \G+i\subset\L_n}\sum_{\D\ni i:\, \D\not\subset\G+i}\sum_{\xi_\D}\int\nu(d\eta)(\frac{1}{|\D|}-\frac{1}{|\D\cap\L_n|})c_\D(\eta,\xi_\D)\log\frac{\mu(\eta_{\D\cap\L_n}|\eta_{\L_n\setminus\D})}{\mu(\xi_{\D\cap\L_n}|\eta_{\L_n\setminus\D})}\cr
&\quad+\sum_{i\in\L_n:\, \G+i\not\subset\L_n}\sum_{\D\ni i}\sum_{\xi_\D}\int\nu(d\eta)(\frac{1}{|\D|}-\frac{1}{|\D\cap\L_n|})c_\D(\eta,\xi_\D)\log\frac{\mu(\eta_{\D\cap\L_n}|\eta_{\L_n\setminus\D})}{\mu(\xi_{\D\cap\L_n}|\eta_{\L_n\setminus\D})}\cr
&=: I+II+III+IV+V,
\end{split}
\end{equation*}
where we used that $\D\subset\L_n$ in $II$.
Now, using once again Lemma~\ref{Lem_Condition}, the terms $I$ and $IV$ can be bounded from above, by 
\begin{equation*}
\begin{split}
|I|\le |\L_n|2\log \frac 1 \d \sum_{\D\ni o:\, \D\not\subset\G}c_\D\qquad\text{ and }\qquad |IV|\le |\L_n|\log \frac 1 \d \sum_{\D\ni o:\, \D\not\subset\G}c_\D
\end{split}
\end{equation*}
where we also applied translation invariance. Since, by well-definedness, $\sum_{\D\ni o}c_\D<\infty$, we can pick $\G\supset\G(\e)$ sufficiently large, such that $|I|+|IV|<\e|\L_n|$.

Further, the boundary terms $III$ and $V$ can be bounded from above, using again translation invariance, by 
\begin{equation*}
\begin{split}
|III|&\le \#\{i\in\L_n:\, \G+i\not\subset\L_n\}2\log \frac 1 \d \sum_{\D\ni o}c_\D\quad\text{ and }\cr
|V|&\le\#\{i\in\L_n:\, \G+i\not\subset\L_n\}\log \frac 1 \d \sum_{\D\ni o}c_\D.
\end{split}
\end{equation*}
For given $\G$, for sufficiently large $n$, we thus have $|III|+|V|<\e|\L_n|$ since $\#\{i\in\L_n:\, \G+i\not\subset\L_n\}/|\L_n|$ tends to zero as $n$ tends to infinity.

In order to bound the remaining bulk term $II$, note that, for any fixed $\D\Subset\Z^d$ with $\D\subset\L_n$, we can estimate using the DLR equation and consistency
\begin{equation*}\label{EntropyChangeAA}
\begin{split}
&\frac{\g_\D(\eta_\D|\eta_{\D^c})\mu(\xi_{\D}|\eta_{\L_n\setminus\D})}{\g_\D(\xi_\D|\eta_{\D^c})\mu(\eta_{\D}|\eta_{\L_n\setminus\D})}
=\frac{\g_\D(\eta_\D|\eta_{\D^c})\mu\big(\g_{\L_n}(\xi_{\D}\eta_{\L_n\setminus\D}|\cdot)\big)}{\g_\D(\xi_\D|\eta_{\D^c})\mu\big(\g_{\L_n}(\eta_{\D}\eta_{\L_n\setminus\D}|\cdot)\big)}\cr
&=\frac{\mu\big(\frac{\g_\D(\xi_{\D}|\eta_{\L_n\setminus\D}\s_{\L_n^c})}{\g_\D(\xi_\D|\eta_{\D^c})}\g_{\L_n}(\eta_{\L_n\setminus\D}|\s_{\L_n^c})\big)}{\mu\big(\frac{\g_\D(\eta_{\D}|\eta_{\L_n\setminus\D}\tilde\s_{\L_n^c})}{\g_\D(\eta_\D|\eta_{\D^c})}\g_{\L_n}(\eta_{\L_n\setminus\D}|\tilde\s_{\L_n^c})\big)}
\le\frac{\sup_{\xi,\eta,\s}\frac{\g_\D(\xi_{\D}|\eta_{\L_n\setminus\D}\s_{\L_n^c})}{\g_\D(\xi_\D|\eta_{\D^c})}}{\inf_{\eta,\s}\frac{\g_\D(\eta_{\D}|\eta_{\L_n\setminus\D}\s_{\L^c})}{\g_\D(\eta_\D|\eta_{\D^c})}}.
\end{split}
\end{equation*}
Hence, using arguments as in the proof of Lemma~\ref{Lem_Condition}, we have $\inf_\eta\g_\D(\eta_\D|\eta_{\D^c})\ge\d^{|\D|}$ and thus
\begin{equation}\label{BulkEstimate}
\begin{split}
1-\d^{-|\D|}|\g_\D(\xi_{\D}|\eta_{\L_n\setminus\D}\s_{\L_n^c})&-\g_\D(\xi_\D|\eta_{\D^c})|\le\frac{\g_\D(\xi_{\D}|\eta_{\L_n\setminus\D}\s_{\L_n^c})}{\g_\D(\xi_\D|\eta_{\D^c})}\cr
&\le 1+\d^{-|\D|}|\g_\D(\xi_{\D}|\eta_{\L_n\setminus\D}\s_{\L_n^c})-\g_\D(\xi_\D|\eta_{\D^c})|
\end{split}
\end{equation}
where l.h.s.~and the r.h.s.~tend to one by the quasilocality assumption on the specification uniformly in the configurations. Thus, for the bulk term $II$, since we are dealing with only finitely many sets $\D\subset\G$, we can pick sufficiently larger $n$ such that 
\begin{equation*}\label{EntropyChangeBBbb}
\begin{split}
\sum_{\D\ni o,\D\subset\G}\sum_{\xi_\D}\frac{1}{|\D|}\int\nu(d\eta)c(\eta,\xi_{\D})|\log\frac{\g_\D(\eta_\D|\eta_{\D^c})\mu(\xi_{\D}|\eta_{\L_n\setminus\D})}{\g_\D(\xi_\D|\eta_{\D^c})\mu(\eta_{\D}|\eta_{\L_n\setminus\D})}|<\e|\L_n|
\end{split}
\end{equation*}
and hence $|II|<\e|\L_n|$.
This finishes the representation part of the proof. 

\medskip
For the continuity let again be $\G\Subset\Z^d$, then 
\begin{equation*}\label{ReprEner}
\begin{split}
\r_L(\nu,\mu)=&\int\nu(d\eta)\sum_{\D\ni o,\D\subset\G}\sum_{\xi_\D}\frac{1}{|\D|}\int c(\eta,\xi_{\D})\log\frac{\g_\D(\eta_\D|\eta_{\D^c})}{\g_\D(\xi_\D|\eta_{\D^c})}\cr
&+\int\nu(d\eta)\sum_{\D\ni o,\D\not\subset\G}\sum_{\xi_\D}\frac{1}{|\D|}\int c(\eta,\xi_{\D})\log\frac{\g_\D(\eta_\D|\eta_{\D^c})}{\g_\D(\xi_\D|\eta_{\D^c})}=:\r_L^\G(\nu)+\r_L^{\G^c}(\nu)
\end{split}
\end{equation*}
and the map $\nu\mapsto\r_L^\G(\nu)$ is weakly continuous as a finite sum of weakly continuous functions by the continuity of the rates and the quasilocality of the specification. The second summand can be bounded from above and below by 
\begin{equation*}\label{ReprEner}
\begin{split}
-\log\frac{1}{\d}\sum_{\D\ni o,\D\not\subset\G}c_\D\le\r_L^{\G^c}(\nu)\le\log\frac{1}{\d}\sum_{\D\ni o,\D\not\subset\G}c_\D
\end{split}
\end{equation*}
which can be made arbitrarily small since we assumed $\sum_{\D\ni o}c_\D<\infty$. Thus $\r_L(\nu,\mu)$ is continuous as a uniform limit of continuous functions.
\end{proof}

\subsection{Proof of Proposition~\ref{ThmUpperSemi}}\label{Thm27}
In perspective of the representation Lemma~\ref{Lem_Rep}, we see that $g_L^n(\nu)$ can be written as
\begin{equation*}\label{EntropyChange1}
\begin{split}
g_L^n(\nu)=\sum_{i\in\L_n}\sum_{\D\ni i}\frac{1}{|\D\cap\L_n|}\sum_{\xi_{\D}}\int\nu(d\eta)c_\D(\eta, \xi_\D)\log\frac{\nu(\xi_{\D\cap\L_n}|\eta_{\L_n\setminus\D})}{\nu(\eta_{\D\cap\L_n}|\eta_{\L_n\setminus\D})}. 
\end{split}
\end{equation*}
Correspondingly, recall the definition of the approximating entropy loss
\begin{equation*}\label{EntropyChange1}
\begin{split}
\tilde g_L^{n}(\nu)=\sum_{i\in\tilde \L_n}\sum_{\D\ni i}\frac{1}{|\D\cap\L_n|}\sum_{\xi_\D}\int\nu(d\eta)c_\D(\eta,\xi_\D)\log\frac{\nu(\xi_{\D\cap\L_n}|\eta_{\L_n\setminus\D})}{\nu(\eta_{\D\cap\L_n}|\eta_{\L_n\setminus\D})}
\end{split}
\end{equation*}
and note that here we just eliminated a part of the first summation. 
The proof of the proposition is a simple application of Lemma~\ref{Lem_Condition}.
\begin{proof}[Proof of Proposition~\ref{ThmUpperSemi}]
Note that by non-nullness of $\nu$ and Lemma~\ref{Lem_Condition}, 
\begin{equation}\label{Bounds_For_g}
\begin{split}
-|\L_n|\log\frac{1}{\d}\sum_{\D\ni o}c_\D\le g_L^n(\nu)\le|\L_n|\log\frac{1}{\d}\sum_{\D\ni o}c_\D.
\end{split}
\end{equation}
Further, since we assume only finitely many types of transitions, see Condition~(\ref{fini}) in Definition~\ref{Gibbs-Attractor-IPS}, for sufficiently large $n$, all $\D\ni i$ with $i\in\tilde\L_n$ and $c_\D>0$ lie in $\L_n$. 
Now, the error $|g_L^{n}(\nu)-\tilde g_L^{n}(\nu)|$ is of boundary order $o(|\L_n|)$, which is immediate from equation \eqref{Bounds_For_g}, compare also the estimates for $III$ and $V$ in the proof of Proposition~\ref{SpecificEnergyLoss}.
\end{proof}

\subsection{Proof of Proposition~\ref{ApproxThmUpperSemi}}\label{Thm29}
The proof rests on a sequence of lemmas which represent separation of harmless components in terms of continuity and a bulk term argument which guarantees the upper-semicontinuity via an application of Jensen's inequality.

\medskip
For convenience let us write $c_\D^{\xi_\D}(\eta):=c_\D(\eta,\xi_\D)$ and recall also the short-hand notations $c_\D(\eta):=\sum_{\xi_\D\neq\eta_\D} c_\D^{\xi_\D}(\eta)$ and $c_\D:=\sup_{\eta} c_\D(\eta)$. 

\medskip
Note that due to Property~(\ref{fini}) in Definition~\ref{Gibbs-Attractor-IPS}, for sufficiently large $n$, in the definition of $\tilde g^n_L$ we can replace $\D\cap\L_n=\D$ for $\D\cap\tilde\L_n\neq\emptyset$ with $c_\D>0$. We will do that in the sequel. In order to create a term that will resembles a relative entropy, we rewrite $\tilde g_L^{n}(\nu)$ as a sum of two terms
\begin{equation}\label{EntropyChange1zt}
\begin{split}
\tilde g_L^{n}(\nu)=&-\sum_{i\in\tilde \L_n}\sum_{\D\ni i}\frac{1}{|\D|}\sum_{\xi_\D}\int\nu(d\eta)c_\D^{\xi_\D}(\eta)\log\frac{\nu(\eta_{\D}|\eta_{\L_n\setminus\D})q^{|\D|}c_\D^{\xi_\D}(\eta)}{\nu(\xi_{\D}|\eta_{\L_n\setminus\D})c_\D(\eta_{\D^c}\xi_\D)}\cr
&+\sum_{i\in\tilde \L_n}\sum_{\D\ni i}\frac{1}{|\D|}\sum_{\xi_\D}\int\nu(d\eta)c_\D^{\xi_\D}(\eta)\log\frac{q^{|\D|}c_\D^{\xi_\D}(\eta)}{c_\D(\eta_{\D^c}\xi_\D)}
=:s^{n}_L(\nu)+r^{n}_L(\nu)
\end{split}
\end{equation}
where the well-definedness of $r^{n}_L(\nu)$ is guaranteed by the no-trap Condition $4$ in Definition~\ref{Gibbs-Attractor-IPS}. In our first supporting lemma we show that the compensating term $r^n_L(\nu)$ has a limiting density which is also continuous. We postpone the proof to the end of this section. 
\begin{lem}\label{SupLem1}
Under the assumptions of Proposition~\ref{ApproxThmUpperSemi} on $L$, for any $\nu\in\PP_\theta$,
$$\lim_{n\uparrow\infty}|\L_n|^{-1}r^{n}_L(\nu)=r_L(\nu)$$ 
exists and $\nu\mapsto r_L(\nu)$ is weakly continuous on $\PP_\theta$. 
\end{lem}
Thus, using Lemma~\ref{SupLem1}, for the proof of Proposition~\ref{ApproxThmUpperSemi}, it suffices to show existence and upper semicontinuity for the density of $s^{n}_L(\nu)$. For this, consider balls $B_n(i):=\{j\in\Z^d: |i-j|\le n\}$ w.r.t.~the Euclidean norm. We define an approximation of $s^{n}_L(\nu)$ given by
\begin{equation*}\label{HSVersion0}
\begin{split}
f^n_L(\nu):=&\sum_{i\in{\tilde\L_n}}\sum_{\D\ni i}\frac{1}{|\D|q^{|\D|}}\sum_{\xi_\D}\sum_{\eta_{\L_n}}\nu(\xi_{\D}\eta_{\L_n\setminus\D})\tilde c_\D(\eta_{B_{n-1}(i)\setminus\D}\xi_{B_{n-1}(i)\cap\D})\cr
&\hspace{3cm}\times\Psi[\frac{1}{\nu(\xi_{\D}\eta_{\L_n\setminus\D})}\int\nu(d\s)\one_{\eta_{\L_n}}(\s)\frac{q^{|\D|}c_\D^{\xi_\D}(\s)}{c_\D(\s_{\D^c}\xi_\D)}].
\end{split}
\end{equation*}
Here, the truncated rates are defined as 
$$\tilde c_\D^{\xi_\D}(\eta_{B_{n-1}(i)}):=\inf_\s c_\D^{\xi_\D}(\eta_{B_{n-1}(i)}\s_{B_{n-1}(i)^c})$$ and in particular, the rate to flip $\tilde c_\D(\eta_{B_{n-1}(i)}):=\sum_{\xi_\D\neq \eta_\D}\tilde c_\D^{\xi_\D}(\eta_{B_{n-1}(i)})$, for $i\in\tilde\L_n$, depends only on the sites $B_{n-1}(i)\subset\L_n$ inside $\L_n$.   
Further, the classical entropy function $\Psi(u):=-u\log u+u-1$ is non-positive and concave, see for example~\cite{Ge11}. 

\medskip
The proof of Proposition~\ref{ApproxThmUpperSemi} will be finished once we have shown the following two lemmas. 
\begin{lem}\label{SupLem2}
Under the assumptions of Proposition~\ref{ApproxThmUpperSemi} on $L$, for any $\nu\in\PP_\theta$,
$$\lim_{n\uparrow\infty}|\L_n|^{-1}f^{n}_L(\nu)=f_L(\nu)\in[-\infty,\infty)$$ 
exists and $\nu\mapsto f_L(\nu)$ is weakly upper-semicontinuous on $\PP_\theta$.
\end{lem}
We postpone the proof of this and the next lemma to the end of this section. 
\begin{lem}\label{SupLem3}
Under the assumptions of Proposition~\ref{ApproxThmUpperSemi} for $L$ for any $\nu\in\PP_\theta$, we have 
$$\lim_{n\uparrow\infty}|\L_n|^{-1}[s^n_L(\nu)-f^{n}_L(\nu)]=0.$$ 
\end{lem}
Let us give now the proof of Proposition~\ref{ApproxThmUpperSemi}.
\begin{proof}[Proof of Proposition~\ref{ApproxThmUpperSemi}]
For any $\nu\in\PP_\theta$, using the Lemmas~\ref{SupLem1}, \ref{SupLem2} and~\ref{SupLem3} we have 
$$\lim_{n\uparrow\infty}|\L_n|^{-1}\tilde g^{n}_L(\nu)=\lim_{n\uparrow\infty}|\L_n|^{-1}[s^{n}_L(\nu)+r^{n}_L(\nu)-f^{n}_L(\nu)+f^{n}_L(\nu)]=r_L(\nu)+f_L(\nu).$$ 
Further, $\nu\mapsto f_L(\nu)$ is weakly upper-semicontinuous and $\nu\mapsto r_L(\nu)$ is weakly continuous on $\PP_\theta$ which implies the desired result.
\end{proof}

\begin{proof}[Proof of Lemma~\ref{SupLem1}]
We can follow similar but much simpler arguments as used in the proof of Proposition~\ref{SpecificEnergyLoss} and show that 
\begin{equation*}
\begin{split}
\lim_{n\uparrow\infty}\frac{1}{|\L_n|}r^{n}_L(\nu)&=\sum_{\D\ni o}\frac{1}{|\D|}\sum_{\xi_\D}\int\nu(d\eta)c_\D^{\xi_\D}(\eta)\log\frac{q^{|\D|}c_\D^{\xi_\D}(\eta)}{c_\D(\eta_{\D^c}\xi_\D)}
=:r_L(\nu).
\end{split}
\end{equation*}
In particular, the mapping $\nu\mapsto r_L(\nu)$ on $\PP_\theta$ is continuous by the continuity of the rates. 
\end{proof}

\begin{proof}[Proof of Lemma~\ref{SupLem2}]
The main argument in the proof is to upper bound $f^n_L$ by $f^{n-1}_L$ using Jensen's inequality via the concavity of $\Psi$. 

\medskip
 Consider $2^d$ disjoined and congruent subcubes $\G_{n,k}$ of $\L_n$ with total side length $2^{n}-1$ as well as $2^d$ disjoined congruent subcubes $\tilde\G_{n,k}$ of $\tilde\L_n$ with total side length $2^{n}-n-1$. Let the subcubes be centered such that $\tilde\G_{n,k}\subset\G_{n,k}$ for each $k$. Note that $\bigcup_{k\in\{1,\dots,2^d\}}\tilde\G_{n,k}\subsetneq\tilde\L_n$. We have
\begin{equation*}\label{Step1_1}
\begin{split}
f^{n}_L(\nu)&\le\sum_{j=1}^{2^d}\sum_{i\in{\tilde\G_{n,j}}}\sum_{\D\ni i}\frac{1}{|\D|q^{|\D|}}\sum_{\xi_\D}\sum_{\eta_{\L_n}}\nu(\xi_{\D}\eta_{\L_n\setminus\D})\tilde c_\D(\eta_{B_{n-1}(i)\setminus\D}\xi_{B_{n-1}(i)\cap\D})\cr
&\hspace{4cm}\times\Psi[\frac{1}{\nu(\xi_{\D}\eta_{\L_n\setminus\D})}\int\nu(d\s)\one_{\eta_{\L_n}}(\s)\frac{q^{|\D|}c_\D^{\xi_\D}(\s)}{c_\D(\s_{\D^c}\xi_\D)}]\cr
&=\sum_{j=1}^{2^d}\sum_{i\in{\tilde\G_{n,j}}}\sum_{\D\ni i}\frac{1}{|\D|q^{|\D|}}\sum_{\xi_\D}\sum_{\eta_{\G_{n,j}}}\tilde c_\D(\eta_{B_{n-1}(i)\setminus\D}\xi_{B_{n-1}(i)\cap\D})\sum_{\eta_{\L_n\setminus\G_{n,j}}}\nu(\xi_{\D}\eta_{\L_n\setminus\D})\cr
&\hspace{4cm}\times\Psi[\frac{1}{\nu(\xi_{\D}\eta_{\L_n\setminus\D})}\int\nu(d\s)\one_{\eta_{\L_n}}(\s)\frac{q^{|\D|}c_\D^{\xi_\D}(\s)}{c_\D(\s_{\D^c}\xi_\D)}]
\end{split}
\end{equation*}
where the first inequality comes from dropping the terms associated to $i\in\tilde\L_n\sm \bigcup_{k\in\{1,\dots,2^d\}}\tilde\G_{n,k}$, since $\Psi$ is non-positive. The equality in the second line is possible since for $i\in\tilde\G_{n,j}$ we have $B_{n-1}(i)\subset\G_{n,j}$ and thus we could move the truncated rates $\tilde c_\D$ in front of the sum over configurations in $\L_n\setminus\G_{n,j}$. 
Now, note that there exists $m\in\N$ such that $c_\D=0$ if $0\in\D\not\subset B_{m-1}(0)$ by the finite range condition on $L$. For $n\ge m$ from $\D\ni i$ and $i\in\tilde\G_{n,j}$ follows $\D\subset\G_{n,j}$. Thus, for $n\ge m$, by an application of Jensen's inequality w.r.t.~the partial sum over configurations in $\L_n\setminus\G_{n,j}$ to the concave function $\Psi$ we have 
\begin{equation*}\label{Step1_1}
\begin{split}
&\sum_{\eta_{\L_n\setminus\G_{n,j}}}\nu(\xi_{\D}\eta_{\L_n\setminus\D})\Psi[\frac{1}{\nu(\xi_{\D}\eta_{\L_n\setminus\D})}\int\nu(d\s)\one_{\eta_{\L_n}}(\s)\frac{q^{|\D|}c_\D^{\xi_\D}(\s)}{c_\D(\s_{\D^c}\xi_\D)}]\cr
&=\nu(\xi_{\D}\eta_{\G_{n,j}\setminus\D})\sum_{\eta_{\L_n\setminus\G_{n,j}}}\frac{\nu(\xi_{\D}\eta_{\L_n\setminus\D})}{\nu(\xi_{\D}\eta_{\G_{n,j}\setminus\D})}\Psi[\frac{1}{\nu(\xi_{\D}\eta_{\L_n\setminus\D})}\int\nu(d\s)\one_{\eta_{\L_n}}(\s)\frac{q^{|\D|}c_\D^{\xi_\D}(\s)}{c_\D(\s_{\D^c}\xi_\D)}]\cr
&\le \nu(\xi_{\D}\eta_{\G_{n,j}\setminus\D})\Psi[\sum_{\eta_{\L_n\setminus\G_{n,j}}}\frac{1}{\nu(\xi_{\D}\eta_{\G_{n,j}\setminus\D})}\int\nu(d\s)\one_{\eta_{\L_n}}(\s)\frac{q^{|\D|}c_\D^{\xi_\D}(\s)}{c_\D(\s_{\D^c}\xi_\D)}]\cr
&= \nu(\xi_{\D}\eta_{\G_{n,j}\setminus\D})\Psi[\frac{1}{\nu(\xi_{\D}\eta_{\G_{n,j}\setminus\D})}\int\nu(d\s)\one_{\eta_{\G_{n,j}}}(\s)\frac{q^{|\D|}c_\D^{\xi_\D}(\s)}{c_\D(\s_{\D^c}\xi_\D)}]
\end{split}
\end{equation*}
which implies that
\begin{equation*}\label{Step1_2}
\begin{split}
f^{n}_L(\nu)\le&\sum_{j=1}^{2^d}\sum_{i\in{\tilde\G_{n,j}}}\sum_{\D\ni i}\frac{1}{|\D|q^{|\D|}}\sum_{\xi_\D}\sum_{\eta_{\G_{n,j}}}\tilde c_\D(\eta_{B_{n-1}(i)\setminus\D}\xi_{B_{n-1}(i)\cap\D})\nu(\xi_{\D}\eta_{\G_{n,j}\setminus\D})\cr
&\times\Psi[\frac{1}{\nu(\xi_{\D}\eta_{\G_{n,j}\setminus\D})}\int\nu(d\s)\one_{\eta_{\G_{n,j}}}(\s)\frac{q^{|\D|}c_\D^{\xi_\D}(\s)}{c_\D(\s_{\D^c}\xi_\D)}]\le 2^d f^{n-1}_L(\nu)
\end{split}
\end{equation*}
where we used translation invariance of $\nu$ and the rates.
Notice that in the last inequality we used that truncating the rates $\tilde c$ over smaller volumes only decreases the rates, which gives the upper bound by non-positivity of $\Psi$.
To compensate for the different volumes define $G(n):=\prod_{l=n}^\infty\frac{(2^{l+2}-2)^d}{(2^{l+2}-1)^d}$ which tends to one for $n\uparrow\infty$, then 
\begin{equation*}\label{Step1_3}
\begin{split}
\frac{G(n)}{(2^{n+1}-1)^d}f^{n}_L(\nu)
\end{split}
\end{equation*}
is non-increasing in $n$ since $f^{n}_L(\nu)\le2^d f^{n-1}_L(\nu)$ and thus 
\begin{equation*}\label{Step1_3}
\begin{split}
\lim_{n\uparrow\infty}\frac{G(n)}{(2^{n+1}-1)^d}f^{n}_L(\nu)=f_L(\nu)\ge-\infty
\end{split}
\end{equation*}
exists. This also implies $\lim_{n\uparrow\infty}|\L_n|^{-1}f^{n}_L(\nu)=f_L(\nu)$. Since $\nu\mapsto f^n_L(\nu)$ is continuous for every $n$, we have that $\nu\mapsto f_L(\nu)$ is upper semicontinuous.
\end{proof}

\begin{proof}[Proof of Lemma~\ref{SupLem3}]
Recall that $\Psi(u)=-u\log u+u-1$. Let us start by decomposing $\Psi$ in $f^n_L(\nu)$ for sufficiently large $n$. We have
\begin{equation}\label{EntropyChangeFullBB}
\begin{split}
f^n_L(\nu)&=\sum_{i\in{\tilde\L_n}}\sum_{\D\ni i}\frac{1}{|\D|q^{|\D|}}\sum_{\xi_\D}\sum_{\eta_{\L_n}}\nu(\xi_{\D}\eta_{\L_n\setminus\D})\tilde c_\D(\eta_{B_{n-1}(i)\setminus\D}\xi_{\D})\cr
&\hspace{3cm}\times\Psi[\frac{1}{\nu(\xi_{\D}\eta_{\L_n\setminus\D})}\int\nu(d\s)\one_{\eta_{\L_n}}(\s)\frac{q^{|\D|}c_\D^{\xi_\D}(\s)}{c_\D(\s_{\D^c}\xi_\D)}]\cr
&=-\sum_{i\in{\tilde\L_n}}\sum_{\D\ni i}\frac{1}{|\D|}\sum_{\xi_\D}\sum_{\eta_{\L_n}}\int\nu(d\s)\one_{\eta_{\L_n}}(\s)\frac{\tilde c_\D(\eta_{B_{n-1}(i)\setminus\D}\xi_{\D})c_\D^{\xi_\D}(\s)}{c_\D(\s_{\D^c}\xi_\D)}\cr
&\hspace{3.5cm}\times\log[\frac{1}{\nu(\xi_{\D}\eta_{\L_n\setminus\D})}\int\nu(d\s)\one_{\eta_{\L_n}}(\s)\frac{q^{|\D|}c_\D^{\xi_\D}(\s)}{c_\D(\s_{\D^c}\xi_\D)}]\cr
&\hspace{0.5cm}+\sum_{i\in{\tilde\L_n}}\sum_{\D\ni i}\frac{1}{|\D|}\sum_{\xi_\D}\int\nu(d\s)\frac{\tilde c_\D(\s_{B_{n-1}(i)\setminus\D}\xi_{\D})c_\D^{\xi_\D}(\s)}{c_\D(\s_{\D^c}\xi_\D)}\cr
&\hspace{0.5cm}-\sum_{i\in{\tilde\L_n}}\sum_{\D\ni i}\frac{1}{|\D|q^{|\D|}}\sum_{\xi_\D}\sum_{\eta_{\L_n}}\nu(\xi_{\D}\eta_{\L_n\setminus\D})\tilde c_\D(\eta_{B_{n-1}(i)\setminus\D}\xi_{\D})\cr
&=:-I+II-III.
\end{split}
\end{equation}
We claim that $II-III$ is of boundary order. Indeed, $II-III$ can be equivalently written as 
\begin{equation}\label{EntropyChangeFullCC}
\begin{split}
\sum_{i\in{\tilde\L_n}}&\sum_{\D\ni i}\frac{1}{|\D|}\Big[\sum_{\xi_\D}\int\nu(d\s)\tfrac{\tilde c_\D(\s_{B_{n-1}(i)\setminus\D}\xi_{\D})c_\D^{\xi_\D}(\s)}{c_\D(\s_{\D^c}\xi_\D)}-\sum_{\eta_{\L_n}}\nu(\eta_{\L_n})\tilde c_\D(\eta_{B_{n-1}(i)})\Big]\cr
&=\sum_{i\in{\tilde\L_n}}\sum_{\D\ni i}\frac{1}{|\D|}\int\nu(d\eta)\Big[\sum_{\xi_\D}\frac{\tilde c_\D(\eta_{B_{n-1}(i)\setminus\D}\xi_{\D})c_\D^{\xi_\D}(\eta)}{c_\D(\eta_{\D^c}\xi_\D)}-\tilde c_\D(\eta_{B_{n-1}(i)})\Big]
\end{split}
\end{equation}
and it suffices to show that the term in square brackets tends to zero as $n\uparrow\infty$ uniformly in $\eta$ and $i$. But this is the case, indeed if $c_\D=0$ by the definition $L$, $\D$ is not included in the summation and hence there is nothing to show. If $c_\D>0$ with $\D\ni i$ we have for all $\eta$ with $c_\D(\eta)>0$, 
\begin{equation*}\label{EntropyChangeCCC}
\begin{split}
&|\sum_{\xi_\D}c_\D^{\xi_\D}(\eta)\frac{\tilde c_\D(\eta_{B_{n-1}(i)\setminus\D}\xi_{\D})}{c_\D(\eta_{\D^c}\xi_\D)}-\tilde c_\D(\eta_{B_{n-1}(i)})|\cr
&\le|\sum_{\xi_\D}c_\D^{\xi_\D}(\eta)\frac{\tilde c_\D(\eta_{B_{n-1}(i)\setminus\D}\xi_{\D})-c_\D(\eta_{\D^c}\xi_\D)}{c_\D(\eta_{\D^c}\xi_\D)}|+\sup_\eta|c_\D(\eta)-\tilde c_\D(\eta_{B_{n-1}(i)})|\cr
&\le\sup_\eta|c_\D(\eta)-\tilde c_\D(\eta_{B_{n-1}(i)})|\sum_{\xi_\D}\frac{c_\D^{\xi_\D}(\eta)}{c_\D(\eta_{\D^c}\xi_\D)}+\sup_\eta|c_\D(\eta)-\tilde c_\D(\eta_{B_{n-1}(i)})|
\end{split}
\end{equation*}
where
\begin{equation*}\label{EntropyChangeCCC}
\begin{split}
\sum_{\xi_\D}\frac{c_\D^{\xi_\D}(\eta)}{c_\D(\eta_{\D^c}\xi_\D)}&=\sum_{\xi_\D:c_\D^{\xi_\D}(\eta)>0}\frac{c_\D^{\xi_\D}(\eta)}{c_\D(\eta_{\D^c}\xi_\D)}\le\frac{c_\D(\eta)}{\min_{\xi_\D:c_\D^{\xi_\D}(\eta)>0}c_\D(\eta_{\D^c}\xi_\D)}\cr
&\le\sup_{\eta:c_\D(\eta)>0}\frac{c_\D(\eta)}{\min_{\xi_\D:c_\D^{\xi_\D}(\eta)>0}c_\D(\eta_{\D^c}\xi_\D)}\le\frac{c_\D}{\inf_{\eta,\xi_\D:c_\D^{\xi_\D}(\eta)>0}c_\D(\eta_{\D^c}\xi_\D)}
\end{split}
\end{equation*}
which is finite by Condition~(\ref{mini}) in Definition~\ref{Gibbs-Attractor-IPS}. Hence, by the uniform continuity of the rates, the density limit of \eqref{EntropyChangeFullCC} tends to zero as $n$ tends to infinity. 

\medskip
It remains to show convergence to zero for the density of $s^{n}_L(\nu)+I$ as $n$ tends to infinity. For this, it suffices to show that, for all $i\in\tilde\L_n$, we have that
\begin{equation*}\label{EntropyChangeD}
\begin{split}
\sum_{\D\ni i}\frac{1}{|\D|}\sum_{\xi_\D}\sum_{\eta_{\L_n}}&\nu(\eta_{\L_n})\Big[\int\nu(d\s|\eta_{\L_n})\frac{\tilde c_\D(\eta_{B_{n-1}(i)\setminus\D}\xi_{\D})c_\D^{\xi_\D}(\s_{\L_n^c}\eta_{\L_n})}{c_\D(\s_{\L_n^c}\eta_{\L_n\sm\D}\xi_\D)}\cr
&\hspace{2cm}\times\log[\int\nu(d\tilde\s|\eta_{\L_n})\frac{\nu(\eta_{\L_n})q^{|\D|}c_\D^{\xi_\D}(\tilde\s_{\L_n^c}\eta_{\L_n})}{\nu(\xi_{\D}\eta_{\L_n\setminus\D})c_\D(\tilde\s_{\L_n^c}\eta_{\L_n\sm\D}\xi_\D)}]\cr
&-\int\nu(d\s|\eta_{\L_n})c_\D^{\xi_\D}(\s_{\L_n^c}\eta_{\L_n})\log\frac{\nu(\eta_{\L_n})q^{|\D|}c_\D^{\xi_\D}(\s_{\L_n^c}\eta_{\L_n})}{\nu(\xi_{\D}\eta_{\L_n\setminus\D})c_\D(\s_{\L_n^c}\eta_{\L_n\sm\D}\xi_\D)}\Big]
\end{split}
\end{equation*}
tends to zero as $n\uparrow\infty$. Adding and subtracting the mixed term 
$$\int\nu(d\s|\eta_{\L_n})c_\D^{\xi_\D}(\s_{\L_n^c}\eta_{\L_n})\log[\int\nu(d\tilde\s|\eta_{\L_n})\frac{\nu(\eta_{\L_n})q^{|\D|}c_\D^{\xi_\D}(\tilde\s_{\L_n^c}\eta_{\L_n})}{\nu(\xi_{\D}\eta_{\L_n\setminus\D})c_\D(\tilde\s_{\L_n^c}\eta_{\L_n\sm\D}\xi_\D)}]$$
we first show boundary order of 
\begin{equation*}\label{EntropyChangeD}
\begin{split}
\sum_{\D\ni i}\frac{1}{|\D|}\sum_{\xi_\D}\sum_{\eta_{\L_n}}\nu(\eta_{\L_n})\int\nu(d\s|\eta_{\L_n})c_\D^{\xi_\D}(\s_{\L_n^c}\eta_{\L_n})\log[\int\nu(d\tilde\s|\eta_{\L_n})\tfrac{c_\D^{\xi_\D}(\tilde\s_{\L_n^c}\eta_{\L_n})c_\D(\s_{\L_n^c}\eta_{\L_n\sm\D}\xi_\D)}{c_\D^{\xi_\D}(\s_{\L_n^c}\eta_{\L_n})c_\D(\tilde\s_{\L_n^c}\eta_{\L_n\setminus\D}\xi_\D)}].
\end{split}
\end{equation*}
Define the strictly positive minimal transition rate guaranteed by Condition~(\ref{mini}) in Definition~\ref{Gibbs-Attractor-IPS} as $c_\D^{\text{min}}$, 
then for the upper bound we have
\begin{equation}\label{EntropyChangeFullDD}
\begin{split}
&\frac{c_\D^{\xi_\D}(\tilde\s_{\L_n^c}\eta_{\L_n})}{c_\D^{\xi_\D}(\s_{\L_n^c}\eta_{\L_n})}\le1+\frac{1}{c_\D^{\text{min}}}\sup_{\eta,\xi,\s}|c_\D^{\xi_\D}(\eta)-c_\D^{\xi_\D}(\eta_{\L_n}\s_{\L_n^c})| \hspace{1cm}\text{ and}\cr
&\frac{c_\D(\s_{\L_n^c}\eta_{\L_n\sm\D}\xi_\D)}{c_\D(\tilde\s_{\L_n^c}\eta_{\L_n\setminus\D}\xi_\D)}\le1 +\frac{1}{c_\D^{\text{min}}}\sup_{\eta,\s}|c_\D(\eta)-c_\D(\eta_{\L_n}\s_{\L_n^c})|
\end{split}
\end{equation}
and similar from below for the lower bound. 
This yields the boundary order by the uniform continuity of the rates.  

\medskip
Secondly, we show boundary order of the second mixed term 
\begin{equation}\label{EntropyChangeJ}
\begin{split}
&\sum_{\D\ni i}\frac{1}{|\D|}\sum_{\xi_\D}\sum_{\eta_{\L_n}}\nu(\eta_{\L_n})\int\nu(d\s|\eta_{\L_n})\Big[\tfrac{\tilde c_\D(\eta_{B_{n-1}(i)\setminus\D}\xi_{\D})c_\D^{\xi_\D}(\s_{\L_n^c}\eta_{\L_n})}{c_\D(\s_{\L_n^c}\eta_{\L_n\sm\D}\xi_\D)}-c_\D^{\xi_\D}(\s_{\L_n^c}\eta_{\L_n})\Big]\cr
&\hspace{3cm}\times\log[\frac{\nu(\eta_{\L_n})}{\nu(\xi_{\D}\eta_{\L_n\setminus\D})}\int\nu(d\tilde\s|\eta_{\L_n})\frac{q^{|\D|}c_\D^{\xi_\D}(\tilde\s_{\L_n^c}\eta_{\L_n})}{c_\D(\tilde\s_{\L_n^c}\eta_{\L_n\sm\D}\xi_\D)}].
\end{split}
\end{equation}
From the argument presented in equation~\eqref{EntropyChangeFullDD} we already see that the term in square brackets in~\eqref{EntropyChangeJ} 
\begin{equation*}\label{EntropyChangeD}
\begin{split}
A(\xi_\D,\eta_{\L_n}):=\int\nu(d\s|\eta_{\L_n})\Big[\tfrac{\tilde c_\D(\eta_{B_{n-1}(i)\setminus\D}\xi_{\D})c_\D^{\xi_\D}(\s_{\L_n^c}\eta_{\L_n})}{c_\D(\s_{\L_n^c}\eta_{\L_n\sm\D}\xi_\D)}-c_\D^{\xi_\D}(\s_{\L_n^c}\eta_{\L_n})\Big]
\end{split}
\end{equation*}
should tend to zero as $n$ tends to infinity uniformly in the configurations by the uniform continuity of the rates. 
Unfortunately, neither the logarithmic term nor $A(\xi_\D,\eta_{\L_n})$ have a fixed sign, which makes a direct estimate difficult. The way out is to split the sums over terms with fixed sign and use two different estimates. One will be $\log x\le x$, the other one is more complicated. Here, we replace the logarithmic term by the function $\Psi$. This has the advantage that the logarithm is then replaced by terms with fixed signs. Reintroducing the function $f^n_L$ will allow us to bound the error against $f^n_L$ itself, which will lead to the following statement. 
If $\lim_{n\uparrow\infty}|\L_n|^{-1}f^n_L(\nu)>-\infty$, then the error tends to zero as $n\uparrow\infty$. In case $\lim_{n\uparrow\infty}|\L_n|^{-1}f^n_L(\nu)=-\infty$ then also $\lim_{n\uparrow\infty}|\L_n|^{-1}s^{n}_L(\nu)=-\infty$. 

\medskip
To do this, let us split the sum $\sum_{\xi_\D,\eta_{\L_n}}$ in \eqref{EntropyChangeJ} into a first sum $\sum_{\xi_\D,\eta_{\L_n}:\, A(\xi_\D,\eta_{\L_n})\ge 0}$ and a second sum $\sum_{\xi_\D,\eta_{\L_n}:\, A(\xi_\D,\eta_{\L_n})< 0}$. The part of~\eqref{EntropyChangeJ} under the second sum can be bounded from above by 
\begin{equation}\label{EntropyChangeHH}
\begin{split}
&\sum_{\D\ni i}\frac{1}{|\D|}\sum_{\xi_\D}\sum_{\eta_{\L_n}}\nu(\xi_\D\eta_{\L_n\sm\D})|A(\xi_\D,\eta_{\L_n})|[\int\nu(d\tilde\s|\eta_{\L_n})\frac{q^{|\D|}c_\D^{\xi_\D}(\tilde\s_{\L_n^c}\eta_{\L_n})}{c_\D(\tilde\s_{\L_n^c}\eta_{\L_n\sm\D}\xi_\D)}]^{-1}
\end{split}
\end{equation}
using $\log x\le x$. The crucial point here is that in this bound, the appearances of $\nu(\eta_{\L_n})$ have cancelled.
Since we assumed the condition of a minimal transition rate, the term $[\int\nu(d\tilde\s|\eta_{\L_n})\frac{c_\D^{\xi_\D}(\tilde\s_{\L_n^c}\eta_{\L_n})}{c_\D(\tilde\s_{\L_n^c}\eta_{\L_n\sm\D}\xi_\D)}]^{-1}$ can be uniformly bounded by $c_\D/c_\D^{\text{min}}$. Then, for sufficiently large $n$, we can bound \eqref{EntropyChangeHH} from above by 
\begin{equation*}\label{EntropyChangeF}
\begin{split}
&\e\sum_{\D\ni o}\frac{c_\D}{c_\D^{\text{min}}},
\end{split}
\end{equation*}
where we also used that $\#\{\eta_\D\}=q^{|\D|}$ and thus there is another cancellation. 
We can use the same arguments to bound the summands $\sum_{\xi_\D,\eta_{\L_n}:\, A(\xi_\D,\eta_{\L_n})\ge 0}$ in~\eqref{EntropyChangeJ} from below by $-\e\sum_{\D\ni o}c_\D/c_\D^{\text{min}}$ and hence they become arbitrarily small as $n\uparrow\infty$.

Now, to bound the summand $\sum_{\xi_\D,\eta_{\L_n}:\, A(\xi_\D,\eta_{\L_n})\ge 0}$ in~\eqref{EntropyChangeJ} from above, using the simple bound $\log x\le x$, the crucial cancellation of terms $\nu(\eta_{\L_n})$ is not available and instead terms involving $\nu$ would remain in the denominator, creating the need to involve non-nullness. In order to circumvent this issue, recall that for $u>0$
$$\log u=1-u^{-1}-u^{-1}\Psi(u)\le 1+u^{-1}-u^{-1}\Psi(u)$$
since only the term $-u^{-1}$ is negative. Using the property of a minimal transition rate, we can even bound $A(\xi_\D,\eta_{\L_n})$ agains its first summand, i.e., for sufficiently large $n$
\begin{equation*}\label{EntropyChangeD}
\begin{split}
|A(\xi_\D,\eta_{\L_n})|\le \e \int\nu(d\s|\eta_{\L_n})\frac{\tilde c_\D(\eta_{B_{n-1}(i)\setminus\D}\xi_{\D})c_\D^{\xi_\D}(\s_{\L_n^c}\eta_{\L_n})}{c_\D(\s_{\L_n^c}\eta_{\L_n\sm\D}\xi_\D)}.
\end{split}
\end{equation*}
Thus, the part of~\eqref{EntropyChangeJ} under this sum can be bounded from above by 
\begin{equation}\label{EntropyChangeH}
\begin{split}
&\e \sum_{\D\ni i}\frac{1}{|\D|}\sum_{\xi_\D}\sum_{\eta_{\L_n}}\nu(\xi_{\D}\eta_{\L_n\setminus\D})\int\nu(d\s|\eta_{\L_n})\frac{\tilde c_\D(\eta_{B_{n-1}(i)\setminus\D}\xi_{\D})c_\D^{\xi_\D}(\s_{\L_n^c}\eta_{\L_n})}{c_\D(\s_{\L_n^c}\eta_{\L_n\sm\D}\xi_\D)}\cr
&\hspace{1cm}\times\Big[\log[\int\nu(d\tilde\s|\eta_{\L_n})\tfrac{q^{|\D|}c_\D^{\xi_\D}(\tilde\s_{\L_n^c}\eta_{\L_n})}{c_\D(\tilde\s_{\L_n^c}\eta_{\L_n\sm\D}\xi_\D)}]+2 [\int\nu(d\tilde\s|\eta_{\L_n})\tfrac{q^{|\D|}c_\D^{\xi_\D}(\tilde\s_{\L_n^c}\eta_{\L_n})}{c_\D(\tilde\s_{\L_n^c}\eta_{\L_n\sm\D}\xi_\D)}]^{-1}\Big].
\end{split}
\end{equation}
But adding and subtracting boundary-order terms of the form $II-III$, this expression is equal to $\e f_L^n(\nu)+o(|\L_n|)$. We can use the same arguments for the remaining case where we have to bound the summands $\sum_{\xi_\D,\eta_{\L_n}:\, A(\xi_\D,\eta_{\L_n})< 0}$ in~\eqref{EntropyChangeJ} from below. 
This completes the proof.
\end{proof}

\subsection{Proof of Theorem~\ref{GTilde}}
The proof uses the detailed balance equations to relate the rates to the specification of the reversible Gibbs measures configuration wise. 
\begin{proof}[Proof of Theorem~\ref{GTilde}]
Using Theorem~\ref{MainCorApprox}, it suffices to prove that Condition~\ref{approx zero entropy loss condition} holds.  For this it is enough to show that $\tilde g_L(\nu|\mu)\le 0$ and that $\tilde g_L(\nu|\mu)=0$ implies that $\nu\in\GG(\g)$.

\medskip
For this, the first part of the proof is similar to the proof of Proposition~\ref{ApproxThmUpperSemi}. Recall that in \eqref{EntropyChange1zt} we write $\tilde g_L^{n}(\nu)$ as a sum of two terms. To simplify notation let us assume $n$ to be sufficiently large such that $\sum_{\D\ni 0: \D\not\subset B_{n-1}(0)}c_\D=0$.
This can be done without loss of generality since we are interested in the large $n$ limit and $L$ is assumed to have the property that there are only finitely many types of transitions.

\medskip
Since we are now in a reversible setting, it is more convenient to extend $\tilde g_L^{n}(\nu)$ in the following way, where we assume $n$ to be sufficiently large,
\begin{equation*}
\begin{split}
\tilde g_L^{\L_n}(\nu)=&-\sum_{i\in\tilde\L_n}\sum_{\D\ni i}\frac{1}{|\D|}\sum_{\xi_\D}\int\nu(d\eta)c_\D^{\xi_\D}(\eta)\log\frac{\nu(\eta_{\L_n})c_\D^{\xi_\D}(\eta)}{\nu(\xi_{\D}\eta_{\L_n\setminus\D})c_\D^{\eta_{\D}}(\eta_{\D^c}\xi_\D)}\cr
&+\sum_{i\in\tilde\L_n}\sum_{\D\ni i}\frac{1}{|\D|}\sum_{\xi_\D}\int\nu(d\eta)c_\D^{\xi_\D}(\eta)\log\frac{c_\D^{\xi_\D}(\eta)}{c_\D^{\eta_{\D}}(\eta_{\D^c}\xi_\D)}
=:s^{n}_L(\nu)+r^{n}_L(\nu)
\end{split}
\end{equation*}
where, by the continuity of the rates
\begin{equation*}\label{EntropyChange1}
\begin{split}
\lim_{n\uparrow\infty}\frac{1}{|\L_n|}r^{n}_L(\nu)&=\sum_{\D\ni o}\frac{1}{|\D|}\sum_{\xi_\D}\int\nu(d\eta)c_\D^{\xi_\D}(\eta)\log\frac{c_\D^{\xi_\D}(\eta)}{c_\D^{\eta_\D}(\eta_{\D^c}\xi_\D)}=:r_L(\nu).
\end{split}
\end{equation*}
Note that $s^{n}_L(\nu)$ is still well-defined since by the reversibility assumption 
$c_\D^{\xi_\D}(\eta)>0$ implies that  $c_\D^{\eta_{\D}}(\eta_{\D^c}\xi_\D)>0$. Indeed, the reversibility implies that for all $\eta_\L$ and $\xi_\D$ with $\D\subset\L$
\begin{equation}\label{Estimate2dg}
\begin{split}
\int\mu(d\s)\g_\L(\eta_\L|\s_{\L^c})c_\D^{\xi_\D}(\s_{\L^c}\eta_\L)=\int\mu(d\s)\g_\L(\eta_{\L\setminus\D}\xi_\D|\s_{\L^c})c_\D^{\eta_\D}(\s_{\L^c}\eta_{\L\setminus\D}\xi_\D).
\end{split}
\end{equation}
Hence, if $c_\D^{\xi_\D}(\eta)>0$ by the continuity also $c_\D^{\xi_\D}(\s_{\L^c}\eta_\L)>0$ for any $\s$, for a sufficiently large volume $\L$. Further, since the specification is assumed to be non-null, also $c_\D^{\eta_\D}(\s_{\L^c}\eta_{\L\setminus\D}\xi_\D)>0$ and $c_\D^{\eta_\D}(\eta_{\D^c}\xi_\D)>0$ for any $\s$, for the same large volume $\L$. 

\medskip
The reversibility in particular implies that $r_L(\nu)+\r_L(\nu,\mu)=0$, i.e.,
\begin{equation*}\label{ReprEner}
\begin{split}
0=\sum_{\D\ni o}\frac{1}{|\D|}\sum_{\xi_\D}\int\nu(d\eta)c_\D^{\xi_\D}(\eta)\log\frac{c_\D^{\xi_\D}(\eta)\g_\D(\eta_\D|\eta_{\D^c})}{c_\D^{\eta_\D}(\eta_{\D^c}\xi_\D)\g_\D(\xi_\D|\eta_{\D^c})}.
\end{split}
\end{equation*}
This can be seen in the following way. As a consequence of \eqref{Estimate2dg} we have 
\begin{equation*}\label{Estimate2gh}
\begin{split}
\frac{\g_\D(\eta_\D|\eta_{\D^c})c_\D^{\xi_\D}(\eta)}{\g_\D(\xi_\D|\eta_{\D^c})c_\D^{\eta_\D}(\eta_{\D^c}\xi_\D)}=\frac{\int\mu(d\s)\frac{\g_\D(\xi_\D|\eta_{\L\setminus\D}\s_{\L^c})c_\D^{\eta_\D}(\s_{\L^c}\eta_{\L\setminus\D}\xi_\D)}{\g_\D(\xi_\D|\eta_{\D^c})c_\D^{\eta_\D}(\eta_{\D^c}\xi_\D)}}{\int\mu(d\s)\frac{\g_\D(\eta_\D|\eta_{\L\setminus\D}\s_{\L^c})c_\D^{\xi_\D}(\s_{\L^c}\eta_\L)}{\g_\D(\eta_\D|\eta_{\D^c})c_\D^{\xi_\D}(\eta)}}
\end{split}
\end{equation*}
where the r.h.s.~tends to one as $\L$ tends to $\Z^d$ by the continuity and non-nullness assumptions on the rates as well as on the specification. 

\medskip
In other words, in a reversible setting, $\tilde g_L(\nu|\mu)=\lim_{n\uparrow\infty}|\L_n|^{-1}s^n_L(\nu)$. Very similar to the proof of Proposition~\ref{ApproxThmUpperSemi}, one can show, using Jensen's inequality, that in the limit as $n$ tends to infinity, $s^n_L(\nu)$ can be replaced by 
\begin{equation*}\label{HSVersion0}
\begin{split}
f^n_L(\nu):=&\sum_{i\in{\tilde\L_n}}\sum_{\D\ni i}\frac{1}{|\D|}\sum_{\xi_\D}\sum_{\eta_{\L_n}}\nu(\xi_{\D}\eta_{\L_n\setminus\D})\tilde c_\D^{\eta_\D}(\eta_{B_{n-1}(i)\setminus\D}\xi_{\D})\cr
&\hspace{3cm}\times\Psi[\frac{1}{\nu(\xi_{\D}\eta_{\L_n\setminus\D})}\int\nu(d\s)\one_{\eta_{\L_n}}(\s)\frac{c_\D^{\xi_\D}(\s)}{c_\D^{\eta_\D}(\s_{\D^c}\xi_\D)}].
\end{split}
\end{equation*}
Again, $a_nf^n_L(\nu)$, with $a_n>0$ some volume-factor, is a non-increasing sequence of non-positive functions. Since $\Psi\le0$ this in particular implies that $\tilde g_L(\nu|\mu)$ exists and $\tilde g_L(\nu|\mu)\le0$, which is the first property that we wanted to check for Condition~\ref{approx zero entropy loss condition}. 

\medskip
As for the second statement, assume that $\tilde g_L(\nu|\mu)=0$ which then implies that $f^n_L(\nu)=0$ for all sufficiently large $n$.
Consequently, for all ${i\in{\tilde\L_n}}$, ${\D\ni i}$, ${\xi_\D}$ and $\eta_{\L_n}$ we have 
\begin{equation}\label{Zero1}
\begin{split}
\nu(\xi_{\D}\eta_{\L_n\setminus\D})\tilde c_\D^{\eta_\D}(\eta_{B_{n-1}(i)\setminus\D}\xi_{\D})\Psi[\frac{1}{\nu(\xi_{\D}\eta_{\L_n\setminus\D})}\int\nu(d\s)\frac{\one_{\eta_{\L_n}}(\xi)c_\D^{\xi_\D}(\s)}{c_\D^{\eta_\D}(\s_{\D^c}\xi_\D)}]=0.
\end{split}
\end{equation}
Let us assume $\tilde c_\D^{\eta_\D}(\eta_{B_{n-1}(i)\setminus\D}\xi_{\D})>0$ and note, as above, that this implies \linebreak $c_\D^{\eta_\D}(\s_{\L_n^c}\eta_{\L_n\setminus\D}\xi_{\D})>0$ and $c_\D^{\xi_\D}(\s_{\L_n^c}\eta_{\L_n})>0$ for all $\s$, by continuity and reversibility.
Under this assumption, $\nu(\eta_{\L_n})=0$ implies $\nu(\xi_{\D}\eta_{\L_n\setminus\D})=0$ since otherwise
\begin{equation*}\label{EntropyChange1} 
\begin{split}
\nu(\xi_{\D}\eta_{\L_n\setminus\D})\tilde c_\D^{\eta_\D}(\eta_{B_{n-1}(i)\setminus\D}\xi_{\D})\Psi[\frac{1}{\nu(\xi_{\D}\eta_{\L_n\setminus\D})}\int\nu(d\s)\frac{\one_{\eta_{\L_n}}(\s)c_\D^{\s_\D}(\s)}{c_\D^{\eta_\D}(\s_{\D^c}\xi_\D)}]<0.
\end{split}
\end{equation*}
In other words, whenever a jump is possible from a configuration $\eta_{\L_n}$ to a configuration $\xi_{\D}\eta_{\L_n\setminus\D}$, then $\nu(\eta_{\L_n})=0$ implies $\nu(\xi_{\D}\eta_{\L_n\setminus\D})=0$. 
By the condition that $L$ is irreducible this implies that from $\nu(\eta_{\L_n})=0$ it follows that $\nu(\xi_{\tilde\L_n}\eta_{\L_n\setminus\tilde\L_n})=0$ for all $\xi_{\tilde\L_n}$. Further assume $\nu(\eta_{\L_n})=0$ for some $\eta_{\L_n}$. Let $m\ge n$ be such that $\tilde\L_m\supset\L_n$, then it follows that $\nu(\xi_{\L_m\setminus\L_n}\eta_{\L_n})=0$ for all $\xi_{\L_m\setminus\L_n}$. Consequently, $\nu(\xi_{\L_m\setminus\L_n}\xi_{\L_n})=0$ for all $\xi_{\L_n}$ and thus $\nu(\eta_{\L_n})=0$ for all $\eta_{\L_n}$ which is a contradiction. Hence $\nu(\eta_{\L_n})>0$ for all $\eta_{\L_n}$.

\medskip
Finally, let $\eta$ by given with $\tilde c_\D^{\eta_\D}(\eta_{B_{n-1}(i)\setminus\D}\xi_{\D})>0$, then using \eqref{Zero1} and the reversibility \eqref{Estimate2dg}, we have 
\begin{equation*}\label{EntropyChange1} 
\begin{split}
1&=\int\nu(d\s|\eta_{\L_n})\tfrac{c_\D^{\xi_\D}(\eta_{\L_n}\s_{\L_n^c})\nu(\eta_\D|\eta_{\L_n\setminus\D})}{c_\D^{\eta_\D}(\s_{\L_n^c}\eta_{\L_n\setminus\D}\xi_\D)\nu(\xi_{\D}|\eta_{\L_n\setminus\D})}\int\nu(d\s|\eta_{\L_n})\tfrac{\g_\D(\xi_\D|\eta_{\L_n\setminus\D}\s_{\L_n^c})\nu(\eta_\D|\eta_{\L_n\setminus\D})}{\g_\D(\eta_\D|\eta_{\L_n\setminus\D}\s_{\L_n^c})\nu(\xi_{\D}|\eta_{\L_n\setminus\D})}.
\end{split}
\end{equation*}
By martingale convergence, this implies that $\nu$ almost surely
\begin{equation}\label{EntropyChange1lkjk}
\begin{split}
\frac{\g_\D(\eta_\D|\eta_{\D^c})}{\g_\D(\xi_\D|\eta_{\D^c})}=\frac{\nu(\eta_{\D}|\eta_{\D^c})}{\nu(\xi_\D|\eta_{\D^c})}.
\end{split}
\end{equation}
Again by the assumption that $L$ is irreducible, the above equation is true for $\nu$ almost all $\eta$ and $\xi_\D\in\{1,\dots, q\}^\D$. Recall the following general fact:
Let $(a_1, ...a_q)$ and $(b_1, ........, b_q)$ be probability vectors with $\frac{a_l}{a_{k}}=\frac{b_l}{b_{k}}$ for all $k,l\in\{1, ........,q\}$ then
\begin{equation*}\label{Glauber_Reversible_is_Gibbs_7}
\begin{split}
a_l=\frac{a_l}{\sum_{k=1}^q a_k}=\frac{1}{1+\sum_{k\neq l} \frac{a_k}{a_l}}=\frac{1}{1+\sum_{k\neq l} \frac{b_k}{b_l}}=b_l.
\end{split}
\end{equation*}
Hence \eqref{EntropyChange1lkjk} implies $\g_\D(\xi_\D|\eta_{\D^c})=\nu(\xi_\D|\eta_{\D^c})$
for $\nu$ almost all $\eta$ and $\xi_\D\in\{1,\dots, q\}^\D$. But this implies that $\nu$ is a Gibbs measure for the specification $\g$.
\end{proof}

\section{Appendix: Independent dynamics}\label{AP}
The first three subsections are completely elementary, and can serve 
as an illustration  which can be read independently, before or after the 
bulk of the paper. 
The last subsection gives an example going beyond our framework, where 
a non-Gibbsian invariant measure occurs caused by lack of reducibility.  

Let us explain the main ideas for the entropy loss and the energy-entropy decomposition first in a 
single-site situation. 
Denote by $\ell$ the generator of a continuous-time Markov chain on the single-site state space 
$\{1,\dots,q\}$ which is irreducible. It is described by the matrix of jump-rates $\ell_{ij}$ 
and acts on test-functions $f$ via
\begin{equation}\label{Ap_1}
\begin{split}
\ell f(i)=\sum_{j=1}^q \ell_{ij}(f(j)-f(i)).
\end{split}
\end{equation}
By irreducibility there is a unique time-stationary (not necessarily reversible) 
distribution given by the probability vector $\mu=(\mu_i)_{i\in 1,\dots,q}$. The 
semigroup giving the probability to jump in time $t$ from $i$ to $j$ is given by the $ij$-th 
entries of the matrix exponential $P^\ell_t= \exp (t\ell)$.

For any single-site starting measure $\nu$, the \textit{relative entropy loss} w.r.t.~
the time-stationary measure is the time-derivative of the relative entropy
\begin{equation*}
\begin{split}
g_\ell(\nu|\mu)&=\frac{d}{dt}_{|t=0}h(\nu P^\ell_t|\mu)=\r_\ell(\nu,\mu)+g_\ell(\nu)
\end{split}
\end{equation*}
with an energetic term and an entropic term 
\begin{equation}\label{Gee}
\begin{split}
\r_\ell(\nu,\mu)&=\sum_{j\neq i}\nu_j\ell_{j i}\log\frac{\mu_j}{\mu_i}\cr
g_\ell(\nu)&=-\sum_{j\neq i}\nu_j\ell_{j i}\log\frac{\nu_j}{\nu_i}.
\end{split}
\end{equation}

\subsection{Detailed balance}
Let us note that in the special case of a reversible dynamics, 
i.e., if $\mu_j\ell_{j i}=\mu_i \ell_{i j}$, we can write the entropy loss 
$$g_\ell(\nu|\mu)=-\sum_{j\neq i}\nu_j\ell_{ji}\log\frac{\nu_j\ell_{ji}}{\nu_i\ell_{ij}}$$
which can be interpreted as a new relative entropy for measures on a doubled system.

The analogue of this step for the proof of semicontinuity in 
the interacting system in finite-volume approximations is
useful as it allows for simplifications in the treatment of boundary terms. This 
makes the treatment of reversible interacting dynamics easier 
than non-reversible dynamics. 

More precisely, defining $\ell^0_{ji}=\ell_{ji}\one_{j\neq i}$ and the pair of measures $\nu^+(j,i)=\nu_j\ell^0_{ji}/Z$ and $\nu^-(j,i)=\nu_i\ell^0_{ij}/Z$ on 
$\{1,\dots,q\}^2$ with normalization $Z=\sum_{j,i}\nu_j\ell^0_{ji}$, we can write
$$g_\ell(\nu|\mu)=-Z h(\nu^+|\nu^-).$$
Note that $g_\ell(\nu|\mu)$ is non-positive.

\subsection{Decomposing entropy loss for non-reversible dynamics}
It turns out to be useful to involve an 
entropy representation also for the non-reversible dynamics, as it suggests a
way to prove the semicontinuity of $\nu\mapsto g_\ell(\nu|\mu)$ also in 
the full infinite-dimensional problem with interacting dynamics.

Using the measures $\nu^\pm$ as defined above, we write the entropic part 
of the entropy loss as 
\begin{equation}\label{Ap_2x}
\begin{split}
g_\ell(\nu)=-Z h(\nu^+|\nu^-)+\sum_{j\neq i}\nu_j\ell_{j i}\log\frac{\ell_{j i}}{\ell_{i j}}.
\end{split}
\end{equation}
The analogue of this decomposition in the infinite-dimensional setting is very 
useful for us to see semicontinuity. Let us outline the reasons. 
Note first that, under sufficient regularity assumptions on $\mu$ and $\ell$, the continuity of the energy part 
$\r_\ell(\nu,\mu)$ holds, as it is a linear function of $\nu$. Also the infinite-volume analogue of 
the second term on the r.h.s.~of~\eqref{Ap_2x} is linear, and 
again under suitable regularity assumptions on $\ell$ is friendly. 
Finally, $h$ is convex, which is an essential ingredient and helps the actual 
proof of upper-semicontinuity in the infinite-volume problem, via finite-volume 
approximations, see the proof of Propositions~\ref{ThmUpperSemi}, \ref{ApproxThmUpperSemi} and Theorem~\ref{GTilde}.

\subsection{Independent infinite-volume dynamics without reversibility assumption}
\label{fuenfdrei}
Let us check the theory in the simple case of independent dynamics. 
Take the independent sum in infinite volume of the generator above
\begin{equation*}
Lf(\eta)=\sum_{x\in\Z^d}\sum_{j=1}^q \ell_{\eta_x,j}[f(\eta^{x,j})-f(\eta)],
\end{equation*}
where $(\eta^{x,j})_y=\eta_y$ for $y\neq x$ and $(\eta^{x,j})_x=j$. 
The associated infinite-volume semigroup $P^L_t$ factorizes over the matrix exponentials $P^\ell_t$.  
Clearly, under irreducibility of the single-site generator $\ell$, 
for any translation-invariant starting measure $\nu$, there is weak convergence of $\nu_t:=\nu P^L_t$ to 
the unique time-stationary measure obtained by tensoring the single-site time-stationary measures over the sites.

This is seen by looking at the probabilities directly, but of course we find from our Theorem
\ref{MainCor}, noting that for independent irreducible  (in the local state spaces) 
dynamics, non-nullness along the trajectory is guaranteed with the bound $\nu_t(\one_{\eta_0}|\eta_{0^c})\ge \min_{\o_0=1,\dots,q}P^\ell_t(\o_0,\eta_0)$ which converges to $\mu_{\eta_0}$ as time tends to infinity.

\subsection{Lack of reducibility, invariant spaces and fuzzy map and potential lack of Gibbsianness}
Let us drop reducibility in the single-site space. We consider an infinite-volume dynamics with  single-site generator which has two disjoint communicating 
classes $C_1,C_2$ so that $C_1\cup C_2=\{1,\dots,q\}$. On each $C_a$ the single-site Markov chain has a time-stationary distribution which we denote by $\mu^a$. No reversibility is assumed.  
Now, for translation-invariant $\nu$, as time tends to infinity, we have weak convergence to the 
$\nu$-dependent time-stationary measure $\mu$ given by the non-trivial mixture  
\begin{equation}\label{fuzzyrandom}
\begin{split}
\nu_t(d\eta) \rightarrow \mu(d\eta):=\int\nu(d\o)\prod_{x\in \Z^d} \Bigl(  
\sum_{a=1,2}\one_{C_a}(\o_x) \mu^a (d\eta_x)
\Bigr).
\end{split}
\end{equation}
Let us introduce the fuzzy map $T(\o_x):=a$ iff $\o_x\in C_a$, extend it to infinite-volume configurations
by $T(\o)_x:=T(\o_x)$, and write $T\mu=\mu\circ T^{-1}$ for  its action on infinite-volume measures.  
The dynamics preserves fuzzy image measure $T\nu $. Hence, the spaces of measures
\begin{equation*}
\begin{split}
D^T_L(\mu):=\{\nu\in\PP_\theta:\, T\nu=T\mu\}
\end{split}
\end{equation*}
are closed under the dynamics. On these spaces weak convergence to $\mu$ indeed takes place. 
While in the case of the reversible particle-exchange dynamics on the lattice as presented in~\cite{Ge79}, under assumption of particle-exchange irreducibility, the family of invariant spaces is indexed by the 
possible single-site distributions on $\{1,\dots,q\}$, in the present case the family of invariant spaces is indexed by the possible 
fuzzy measures on $\{1,2\}^{\Z^d}$. 

Let us finally mention that there are particular examples for which the limiting measure $\mu$ is non-Gibbsian, even if the starting measure $\nu$ is Gibbsian: 
Consider a $2k$-state Potts model $\nu$ at low temperatures in zero magnetic field. To be specific assume that $\nu$ is 
obtained with fixed boundary condition, say equal to $1$. 
Put $C_1=\{1,\dots,k\}$, $C_2=\{k+1,\dots, 2k\}$.  Consider any independent dynamics with these 
communicating classes, which may not be reversible (assuming $k\geq 3$).  
Then, the measure $T\nu$, which is known as a fuzzy Potts measure, 
is proved to be non-Gibbs at low enough temperatures, see~\cite{Ha03}, for a lack of quasilocality.
Observe, that the conditional probabilities of the time-stationary measure $\mu$ and 
the fuzzy measure $T\nu$ corresponding to  the starting measure $\nu$, are related via 
\begin{equation*}\label{fuzzymakesnG}
\begin{split}
\mu(\one_{\eta_0}|\eta_{V\ba 0})=T\nu(\one_{T(\eta_0)}|T(\eta_{V\ba 0}))\mu^{T(\eta_0)}(\one_{\eta_0}).
\end{split}
\end{equation*}
In this case, also the corresponding time-stationary infinite-volume measure $\mu$ 
in~\eqref{fuzzyrandom} is a non-Gibbsian measure. This goes beyond our original framework and points towards another line of research where (dependent) 
dynamics are considered that have no Gibbsian time-stationary measure, where one may 
hope to make progress by a combination of 
the present methods and those of \cite{KuLeRe04} with new ideas.

\end{document}